\documentclass[final,1p,times]{elsarticle}

\usepackage{amssymb}
 \usepackage{amsthm}
\usepackage{amscd}
\usepackage{amsmath}
\usepackage{amsfonts}
\usepackage{amssymb}
\usepackage{graphicx}
\newtheorem{theorem}{Theorem}

\newtheorem{lemma}[theorem]{Lemma}
\newtheorem{corollary}[theorem]{Corollary}
\usepackage{mathrsfs}
\usepackage{titletoc}

\newcommand{\ra}{\rightarrow}
\newcommand{\p}{\partial}
\newcommand{\f}{\frac}

\newcommand{\be}{\begin{equation}}
\renewcommand{\ra}{\rightarrow}
\newcommand{\ee}{\end{equation}}
\newcommand{\bea}{\begin{eqnarray}}
\newcommand{\eea}{\end{eqnarray}}
\newcommand{\bna}{\begin{eqnarray*}}
\newcommand{\ena}{\end{eqnarray*}}

\renewcommand{\le}{\left}
\newcommand{\ri}{\right}

\journal{***}

\begin{document}

\begin{frontmatter}

\title{Extremal functions for singular Trudinger-Moser inequalities in
 the entire Euclidean space}

\author{Xiaomeng Li$^{1,2}$}
\ead{xmlimath@ruc.edu.cn}

\author{Yunyan Yang$^1$}
 \ead{yunyanyang@ruc.edu.cn}
\address{$^1$Department of Mathematics,
Renmin University of China, Beijing 100872, P. R. China}
\address{$^2$ School of Information, Huaibei Normal University, Huaibei, 235000, P. R. China}

\begin{abstract}
In a previous work (Int. Math. Res. Notices 13 (2010) 2394-2426), Adimurthi-Yang proved a singular Trudinger-Moser inequality in the entire
Euclidean space $\mathbb{R}^N$ $(N\geq 2)$. Precisely, if $0\leq \beta<1$ and $0<\gamma\leq1-\beta$, then there holds for any $\tau>0$,
$$\sup_{u\in W^{1,N}(\mathbb{R}^N),\,\int_{\mathbb{R}^N}(|\nabla u|^N+\tau |u|^N)dx\leq 1}\int_{\mathbb{R}^N}
\f{1}{|x|^{N\beta}}\le(e^{\alpha_N\gamma|u|^{\f{N}{N-1}}}-\sum_{k=0}^{N-2}\f{\alpha_N^k\gamma^k|u|^{\f{kN}{N-1}}}{k!}\ri)dx<\infty,$$
where $\alpha_N=N\omega_{N-1}^{1/(N-1)}$ and $\omega_{N-1}$ is the area of
the unit sphere in $\mathbb{R}^N$.
The above inequality is sharp in the sense that if $\gamma>1-\beta$, all integrals are still finite but the supremum is infinity.
In this paper, we concern extremal functions for these singular inequalities. The regular case $\beta=0$ has been considered by Li-Ruf
(Indiana Univ. Math. J. 57 (2008) 451-480) and Ishiwata (Math. Ann. 351 (2011) 781-804). We shall investigate the singular case $0<\beta<1$ and
prove that for all $\tau>0$, $0<\beta<1$ and $0<\gamma\leq 1-\beta$, extremal functions for the above inequalities exist. The proof is based on
blow-up analysis.
\end{abstract}


\end{frontmatter}

\section{Introduction and main results}
Let $\Omega\subset\mathbb{R}^N$ $(N\geq 2)$ be a bounded smooth domain, $W_0^{1,N}(\Omega)$ be the usual
Sobolev space. Denote $\alpha_N=N\omega_{N-1}^{1/(N-1)}$, where $\omega_{N-1}$ is the area of the
unit sphere in $\mathbb{R}^N$. The famous Trudinger-Moser inequality \cite{24,19,17,22,14} reads
 \be\label{Tru-Moser}\sup_{u\in W_0^{1,N}(\Omega),\,\int_\Omega|\nabla u|^Ndx\leq 1}
 \int_\Omega e^{\alpha |u|^{\f{N}{N-1}}}dx<\infty, \quad
 \forall \alpha\leq \alpha_N.\ee
 This inequality is sharp in the sense that all integrals are still finite when $\alpha>\alpha_N$, but the supremum
 is infinity. It was extended by Cao \cite{Cao}, J. M. do \'O \cite{doo1}, Panda \cite{Panda},
 Ruf \cite{Ruf}, and Li-Ruf \cite{Li-Ruf} to the entire Euclidean space $\mathbb{R}^N$ $(N\geq 2)$. Namely
 \be\label{Tru-Moser-RN}\sup_{u\in W^{1,N}(\mathbb{R}^N),\,\int_{\mathbb{R}^N}(|\nabla u|^N+|u|^N)dx\leq 1}
 \int_{\mathbb{R}^N} \le(e^{\alpha |u|^{\f{N}{N-1}}}-\sum_{k=0}^{N-2}\f{\alpha^k|u|^{\f{Nk}{N-1}}}{k!}\ri)dx<\infty, \quad
 \forall \alpha\leq \alpha_N.\ee
 Recently, several interesting developments of (\ref{Tru-Moser-RN}) has been obtained by J. M. do \'O and M. de Souza \cite{doo2,doo3}.

 Using a rearrangement argument and a change of variables, Adimurthi-Sandeep \cite{A-S} generalized the Trudinger-Moser inequality
 (\ref{Tru-Moser}) to
 a singular version as follows:
 \be\label{Singular}\sup_{u\in W_0^{1,N}(\Omega),\,\int_\Omega |\nabla u|^Ndx\leq 1}
 \int_\Omega \f{e^{\alpha_N \gamma |u|^{\f{N}{N-1}}}}{|x|^{N\beta}}dx<\infty, \quad
 0\leq \beta<1,\,\, 0<\gamma\leq 1-\beta.\ee
 This inequality is also sharp in the sense that all integrals are still finite when $\gamma>1-\beta$, but the supremum
 is infinity. Obviously, if $\beta=0$, then (\ref{Singular}) reduces to (\ref{Tru-Moser}). Later,
 (\ref{Singular}) was extended to the entire $\mathbb{R}^N$ by Adimurthi-Yang \cite{Adi-Yang}.
 Precisely there holds for constants $\tau>0$, $0\leq\beta<1$ and $0<\gamma\leq 1-\beta$,
 \be\label{Singular-RN}\sup_{\int_{\mathbb{R}^N}(|\nabla
  u|^N+\tau|u|^N)dx\leq 1}
 \int_{\mathbb{R}^N}\f{1}{|x|^{N\beta}}\le(e^{\alpha_N\gamma|u|^{N/(N-1)}}-\sum_{k=0}^{N-2}
 \f{(\alpha_N\gamma)^k|u|\,^{kN/(N-1)}}{k!}\ri)dx<\infty.\ee
 Clearly, (\ref{Tru-Moser-RN}) is a special case of (\ref{Singular-RN}). It should be remarked that in \cite{Adi-Yang},
 the proof of (\ref{Singular-RN}) is essentially based on the Young inequality; while in \cite{Li-Ruf}, (\ref{Tru-Moser-RN}) is proved
 via the method of blow-up analysis.  Such kind of singular Trudinger-Moser inequalities are very important in analysis of partial
 differential equations, see for examples \cite{YangJFA,YangJDE,Yangjfa}.

An interesting problem on Trudinger-Moser inequalities is whether or not extremal functions exist.
Existence of extremal functions for the Trudinger-Moser inequality (\ref{Tru-Moser}) was obtained by Carleson-Chang \cite{C-C}
when $\Omega$ is the unit ball, by M. Struwe \cite{Struwe} when $\Omega$ is close to the ball in the sense of measure,
  by M. Flucher and K. Lin
\cite{Flucher,Lin} when $\Omega$ is a general bounded smooth domain, and by Y. Li \cite{Lijpde}
for compact Riemannian surfaces.
For recent developments, we refer the reader to Yang \cite{YangJDE-2015}.
On extremal functions for (\ref{Tru-Moser-RN}), it was proved by Ruf \cite{Ruf} and  Ishiwata \cite{Ishiwata} that if
$N=2$, then there exists some $\epsilon_0>0$ such that for all $\epsilon_0<\alpha\leq 2\pi$, the supremum
$$\sup_{u\in W^{1,2}(\mathbb{R}^2),\,\int_{\mathbb{R}^2}(|\nabla u|^2+u^2)dx\leq 1}\int_{\mathbb{R}^2}(e^{\alpha u^2}-1)dx$$
can be attained by some function $u\in W^{1,2}(\mathbb{R}^2)$ satisfying $\|u\|_{W^{1,2}(\mathbb{R}^2)}\leq 1$.
While for sufficiently small $\alpha>0$, the above supremum can not be attained. If $N\geq 3$, then for any $0\leq\alpha<\alpha_N$,
the supremum in (\ref{Tru-Moser-RN})
can be achieved. While Li-Ruf \cite{Li-Ruf} proved that  when $\alpha=\alpha_N$, extremal function exists for the above supremum.

Our aim is to find extremal functions for the singular Trudinger-Moser inequality (\ref{Singular-RN}) in the case $0<\beta<1$.
Note that the case $\beta=0$ has been studied by Ruf \cite{Ruf}, Ishiwata \cite{Ishiwata} and Li-Ruf \cite{Li-Ruf}. While
these two situations are quite different in analysis.
Throughout this paper, we write for all $\tau\in (0,\infty)$,
\be\label{1tau}\|u\|_{1,\tau}=\le(\int_{\mathbb{R}^N}|\nabla u|^Ndx+\tau\int_{\mathbb{R}^N}|u|^Ndx\ri)^{1/N}.\ee
Obviously $\|\cdot\|_{1,\tau}$ is equivalent to the standard Sobolev norm on $W^{1,N}(\mathbb{R}^N)$.
Define a function $\zeta:\mathbb{N}\times\mathbb{R}\ra \mathbb{R}$ by
\be\label{MF}\zeta(N,s)=e^s-\sum_{k=0}^{N-2}\f{s^k}{k!}=\sum_{k=N-1}^{\infty}\f{s^k}{k!}.\ee

Our main results are the existence of extremal functions for subcritical or critical singular Trudinger-Moser inequality, which can be
stated as the following two theorems respectively.
\begin{theorem}{\bf (Subcritical case)}\label{Theorem 1}
Let $N\geq 2$, $\tau>0$, $\|\cdot\|_{1,\tau}$ and $\zeta:\mathbb{N}\times\mathbb{R}\ra \mathbb{R}$ be defined as in (\ref{1tau}) and
(\ref{MF}) respectively. Then for any $0<\beta<1$ and $0<\epsilon<1-\beta$, the supremum
\be\label{subcrit}\Lambda_{N,\beta,\tau,\epsilon}=\sup_{u\in W^{1,N}(\mathbb{R}^N),\,\|u\|_{1,\tau}\leq 1}
 \int_{\mathbb{R}^N}\f{\zeta(N,\alpha_N(1-\beta-\epsilon)|u|^{\f{N}{N-1}})}{|x|^{N\beta}}dx\ee
 can be attained by some nonnegative decreasing radially symmetric function $u_\epsilon\in C^1(\mathbb{R}^N\setminus
\{0\})\cap C^0(\mathbb{R}^N)\cap W^{1,N}(\mathbb{R}^N)$ with $\|u_\epsilon\|_{1,\tau}=1$.
\end{theorem}

\begin{theorem}{\bf (Critical case)}\label{Theorem 2}
Let $N\geq 2$, $\tau>0$, $\|\cdot\|_{1,\tau}$ and $\zeta:\mathbb{N}\times\mathbb{R}\ra \mathbb{R}$ be defined as in (\ref{1tau}) and
(\ref{MF}) respectively. Then for any $0<\beta<1$, the supremum
\be\label{crit}\Lambda_{N,\beta,\tau}=\sup_{u\in W^{1,N}(\mathbb{R}^N),\,\|u\|_{1,\tau}\leq 1}
 \int_{\mathbb{R}^N}\f{\zeta(N,\alpha_N(1-\beta)|u|^{\f{N}{N-1}})}{|x|^{N\beta}}dx\ee
 can be attained by some nonnegative decreasing radially symmetric function $u^\ast\in C^1(\mathbb{R}^N\setminus
\{0\})\cap C^0(\mathbb{R}^N)\cap W^{1,N}(\mathbb{R}^N)$ with $\|u^\ast\|_{1,\tau}=1$.
\end{theorem}

Trudinger-Moser inequalities involved the norm $\|\cdot\|_{1,\tau}$ was first
introduced by Adimurthi-Yang \cite{Adi-Yang}. This type of inequalities are easy to use in analysis of partial
differential equations with exponential growth. It should be remarked that both the above inequalities
and existence of extremal functions are independent of $\tau$. Let us give the outline of proving Theorems
\ref{Theorem 1} and \ref{Theorem 2}.
The proof of Theorem \ref{Theorem 1}
is based on a direct method of  variation.
By a rearrangement argument, we can take a maximizing sequence $u_j$ satisfying $u_j\geq 0$ and decreasing radially symmetric.
Clearly $u_j\rightharpoonup u_\epsilon$ weakly in $W^{1,N}(\mathbb{R}^N)$ for some $u_\epsilon$.
Since $0<\epsilon<1-\beta$ and $0<\beta<1$, for any $\nu>0$, there exists sufficiently large $R>0$ such that
$$\int_{|x|>R}\f{\zeta(N,\alpha_N(1-\beta-\epsilon)|u_j|^{\f{N}{N-1}})}{|x|^{N\beta}}dx<\nu.$$
 Since $\alpha_N(1-\beta-\epsilon)<\alpha_N(1-\beta)$, we have by the singular Trudinger-Moser inequality
 (\ref{Singular-RN}) that
 $$\lim_{j\ra\infty}\int_{|x|\leq R}\f{\zeta(N,\alpha_N(1-\beta-\epsilon)|u_j|^{\f{N}{N-1}})}{|x|^{N\beta}}dx=
 \int_{|x|\leq R}\f{\zeta(N,\alpha_N(1-\beta-\epsilon)|u_\epsilon|^{\f{N}{N-1}})}{|x|^{N\beta}}dx.$$
Then the conclusion of Theorem \ref{Theorem 1} follows from the above two estimates.

Following Li-Ruf \cite{Li-Ruf} and thereby  following closely Carleson-Chang \cite{C-C}, Ding-Jost-Li-Wang \cite{DJLW} and
Adimurthi-Struwe \cite{Adi-Str},
we prove Theorem \ref{Theorem 2} via the method of blow-up analysis.  Particularly we divide the proof into several steps:

{\it Step 1.} For any $0<\epsilon<1-\beta$, the supremum $\Lambda_{N,\beta,\tau,\epsilon}$
 can be attained by some function $u_\epsilon$ (This is the content of Theorem \ref{Theorem 1} exactly). The Euler-Lagrange
equation of $u_\epsilon$ is semi-linear elliptic when $N=2$, or quasi-linear elliptic when $N\geq 3$;

{\it Step 2.}
Denote $c_\epsilon=u_\epsilon(0)=\max_{\mathbb{R}^N}u_\epsilon$. If $c_\epsilon$ is a bounded sequence, then applying elliptic estimates
to the equation of $u_\epsilon$, we conclude that $u_\epsilon$ converges to a desired extremal function in
$C^1_{\rm loc}(\mathbb{R}^N\setminus\{0\})\cap C^0_{\rm loc}(\mathbb{R}^N)$. If $c_\epsilon\ra+\infty$, then by a delicate analysis
on $u_\epsilon$, we derive
$$\Lambda_{N,\beta,\tau}=\lim_{\epsilon\ra 0}\int_{\mathbb{R}^N}
\f{\zeta(N,\alpha_N(1-\beta-\epsilon)u_\epsilon^{\f{N}{N-1}})}{|x|^{N\beta}}dx\leq \f{1}{1-\beta}\f{\omega_{N-1}}{N}e^{\sum_{k=1}^{N-1}\f{1}{k}+\alpha_N(1-\beta)A_0}.$$
Here $A_0=\lim_{x\ra 0}(G(x)+(N/\alpha_N) \log|x|)$, $G$ is a Green function satisfying
$$-{\rm div}(|\nabla G|^{N-2}\nabla G)+\tau G^{N-1}=\delta_0\quad{\rm in}\quad \mathbb{R}^N,$$
where $\delta_0$ is a Dirac measure centered at $0$.

{\it Step 3.} We construct a sequence of functions $\phi_\epsilon\in W^{1,N}(\mathbb{R}^N)$ satisfying $\|\phi_\epsilon\|_{1,\tau}=1$
and if $\epsilon$ is sufficiently small, then
$$\int_{\mathbb{R}^N}\f{\zeta(N,\alpha_N(1-\beta)\phi_\epsilon^{\f{N}{N-1}})}{|x|^{N\beta}}dx>
\f{1}{1-\beta}\f{\omega_{N-1}}{N}e^{\sum_{k=1}^{N-1}\f{1}{k}+\alpha_N(1-\beta)A_0}.$$

Comparing {\it Steps}  2 and 3, we conclude that $c_\epsilon$ must be bounded and thus the existence of extremal function
follows from elliptic estimates. It should be remarked that in {\it Step 2}, we shall use an estimate of Carleson-Chang \cite{C-C}:
\begin{lemma}\label{C-C-Lemma}
Let $B_1$ be the unit ball in $\mathbb{R}^N$, $v_\epsilon\in W_0^{1,N}(B_1)$ satisfy $\int_{B_1}|\nabla v_\epsilon|^Ndx\leq 1$,
and $v_\epsilon\rightharpoonup 0$
weakly in $W_0^{1,N}(B_1)$. Then
$$
\limsup_{\epsilon\ra 0}\int_{B_1}(e^{\alpha_N |v_\epsilon|^{{N}/{(N-1)}}}-1)dx\leq
\f{\omega_{N-1}}{N}e^{\sum_{k=1}^{N-1}\f{1}{k}}.
$$
\end{lemma}

Before ending this introduction, we mention Csato-Roy \cite{CR}, Iula-Mancini \cite{I-M} and
Yang-Zhu \cite{Yang-Zhu} who studied the same topic in bounded planar domain or compact Riemannian surface.
Throughout this paper, we do {\it not} distinguish sequence and subsequence, the reader can easily
see it from the context. We denote a ball centered at $0$ with radius $r$ by $B_r$,
 $o_\epsilon(1)\ra 0$ as $\epsilon\ra 0$, $o_r(1)\ra 0$ as $r\ra 0$, and
$o_R(1)\ra 0$ as $R\ra\infty$. \\

The remaining part of this paper is devoted to the proof of Theorems \ref{Theorem 1} and \ref{Theorem 2} and organized as follows:
Since the proof is transparent in $\mathbb{R}^2$, we show it in Section \ref{Sec 2}.
In Section \ref{Sec 3}, we prove Theorems \ref{Theorem 1} and \ref{Theorem 2} in $N(\geq 3)$ dimensions.

\section{Two dimensional case}\label{Sec 2}

When $N=2$, extremal functions for subcritical singular Trudinger-Moser inequalities are distributional
  solutions of elliptic partial differential equations of second order. Compared with $N\geq 3$, analysis in two dimensions
   becomes much easier and transparent, so we deal with this case first.

\subsection{Proof of Theorem \ref{Theorem 1}}

We rephrase Theorem \ref{Theorem 1} as below:
\begin{theorem}\label{prop3}
Let $\tau>0$ and $0<\beta<1$ be fixed. Then for any $0<\epsilon<1-\beta$, there exists some nonnegative
decreasing radially symmetric function $u_\epsilon\in C^1(\mathbb{R}^2\setminus\{0\})\cap C^0(\mathbb{R}^2)\cap W^{1,2}(\mathbb{R}^2)$
satisfying $\|u_\epsilon\|_{1,\tau}=1$ and
\be\label{sct}\int_{\mathbb{R}^2}
\f{e^{4\pi(1-\beta-\epsilon)u_\epsilon^2}-1}{|x|^{2\beta}}dx=\Lambda_{2,\beta,\tau,\epsilon}=
\sup_{u\in W^{1,2}(\mathbb{R}^2),\,\|u\|_{1,\tau}\leq 1}\int_{\mathbb{R}^2}
\f{e^{4\pi(1-\beta-\epsilon)u^2}-1}{|x|^{2\beta}}dx.\ee
\end{theorem}

{\it Proof.} Let $\tau>0$, $0<\beta<1$ and $0<\epsilon<1-\beta$ be fixed.
Suppose that $\widetilde{u}$ is the decreasing rearrangement of $|u|$. It is known that
$\int_{\mathbb{R}^2}\widetilde{u}^2dx=\int_{\mathbb{R}^2}{u}^2dx$,
$\int_{\mathbb{R}^2}|\nabla \widetilde{u}|^2dx\leq \int_{\mathbb{R}^2}|\nabla {u}|^2dx$ and
$$\int_{\mathbb{R}^2}\f{e^{4\pi(1-\beta-\epsilon)\widetilde{u}^2}-1}{|x|^{2\beta}}dx\geq\int_{\mathbb{R}^2}
\f{e^{4\pi(1-\beta-\epsilon)u^2}-1}{|x|^{2\beta}}dx.$$
Here we used the Hardy-Littlewood inequality in the last estimate. Therefore
we have
$$\Lambda_{2,\beta,\tau,\epsilon}=\sup_{u\in\mathscr{S}}\int_{\mathbb{R}^2}
\f{e^{4\pi(1-\beta-\epsilon)u^2}-1}{|x|^{2\beta}}dx,$$
where
$\mathscr{S}$ is a set consisting of all nonnegative decreasing radially symmetric functions $u\in W^{1,2}(\mathbb{R}^2)$ with
$\|u\|_{1,\tau}\leq 1$. Take $u_j\in\mathscr{S}$ such that $\int_{\mathbb{R}^2}(e^{4\pi(1-\beta-\epsilon)u_j^2}-1)/|x|^{2\beta}dx
\ra\Lambda_{2,\beta,\tau,\epsilon}$ as $j\ra\infty$. Without loss of generality, we
can assume that there exists some function $u_\epsilon\in W^{1,2}(\mathbb{R}^2)$ such that up to
 a subsequence, as $j\ra\infty$, there holds $u_j\rightharpoonup u_\epsilon$ weakly in $W^{1,2}(\mathbb{R}^2)$, $u_j\ra u_\epsilon$ in
 $L^p_{\rm loc}(\mathbb{R}^2)$ for any $p>0$ and $u_j\ra u_\epsilon$ a.e. in $\mathbb{R}^2$. Hence up to a set of zero measure, $u_\epsilon$ is nonnegative
 decreasing radially symmetric on $\mathbb{R}^2$. Moreover, we have that $\|u_\epsilon\|_{1,\tau}\leq \limsup_{j\ra\infty}\|u_j\|_{1,\tau}\leq 1$.
 Note that $0<\beta<1$ and $0<\epsilon<1-\beta$. Given any $\nu>0$, in view of the Trudinger-Moser inequality (\ref{Tru-Moser-RN}),
 there exists a sufficiently large $r>0$ such that for all $u\in W^{1,2}(\mathbb{R}^2)$ with $\|u\|_{1,\tau}\leq 1$,
 \be\label{small}\f{1}{r^{2\beta}}\int_{|x|> r}(e^{4\pi(1-\beta-\epsilon)u^2}-1)dx\leq
 \f{1}{r^{2\beta}}\int_{\mathbb{R}^2}(e^{4\pi(1-\beta-\epsilon)u^2}-1)dx<\nu.\ee
 Since $u_j\ra u_\epsilon$ in $L^p_{\rm loc}(\mathbb{R}^2)$ for any $p>0$, we have by using the mean value theorem,
 \be\label{le-r}\int_{|x|\leq r}\f{e^{4\pi(1-\beta-\epsilon)u_\epsilon^2}-1}{|x|^{2\beta}}dx=\lim_{j\ra\infty}
 \int_{|x|\leq r}\f{e^{4\pi(1-\beta-\epsilon)u_j^2}-1}{|x|^{2\beta}}dx.\ee
 Combining (\ref{small}) and (\ref{le-r}), we obtain
 $$\int_{\mathbb{R}^2}
\f{e^{4\pi(1-\beta-\epsilon)u_\epsilon^2}-1}{|x|^{2\beta}}dx-\nu\leq\limsup_{j\ra\infty}\int_{\mathbb{R}^2}
\f{e^{4\pi(1-\beta-\epsilon)u_j^2}-1}{|x|^{2\beta}}dx\leq \int_{\mathbb{R}^2}
\f{e^{4\pi(1-\beta-\epsilon)u_\epsilon^2}-1}{|x|^{2\beta}}dx+\nu.$$
Since $\nu$ is arbitrary, there holds
$$\lim_{j\ra\infty}\int_{\mathbb{R}^2}
\f{e^{4\pi(1-\beta-\epsilon)u_j^2}-1}{|x|^{2\beta}}dx=\int_{\mathbb{R}^2}
\f{e^{4\pi(1-\beta-\epsilon)u_\epsilon^2}-1}{|x|^{2\beta}}dx.$$
This leads to (\ref{sct}). Noting that
$$\int_{\mathbb{R}^2}
\f{e^{4\pi(1-\beta-\epsilon)u_\epsilon^2}-1}{|x|^{2\beta}}dx\leq \int_{\mathbb{R}^2}
\f{e^{4\pi(1-\beta-\epsilon)\f{u_\epsilon^2}{\|u_\epsilon\|_{1,\tau}^2}}-1}{|x|^{2\beta}}dx,$$
we get the extremal function $u_\epsilon$, which is nonnegative and decreasing radially symmetric, and satisfies
$\|u_\epsilon\|_{1,\tau}=1$. A straightforward calculation shows that $u_\epsilon$ satisfies the following Euler-Lagrange equation
\be\label{Euler-Lagrange}\le\{\begin{array}{lll}-\Delta u_\epsilon+\tau u_\epsilon=\f{1}{\lambda_\epsilon}\f{u_\epsilon
e^{4\pi(1-\beta-\epsilon) u_\epsilon^2}}{|x|^{2\beta}}\quad{\rm in}\quad \mathbb{R}^2,\\[1.2ex]
u_\epsilon> 0\quad {\rm in}\quad\mathbb{R}^2,\\[1.2ex]
\|u_\epsilon\|_{1,\tau}=1,
\\[1.2ex]
\lambda_\epsilon=\int_{\mathbb{R}^2} |x|^{-2\beta}u_\epsilon^2e^{4\pi(1-\beta-\epsilon) u_\epsilon^2}dx.
\end{array}\ri.\ee
 Applying elliptic estimates to (\ref{Euler-Lagrange}), we have
$u_\epsilon\in C^1(\mathbb{R}^2\setminus\{0\})\cap C^0(\mathbb{R}^2)$.
Here $u_\epsilon>0$ follows from the classical maximum principle and the fact that
$u_\epsilon(0)=\max_{\mathbb{R}^2}u_\epsilon$. This completes the proof of the theorem.
$\hfill\Box$\\

From now on, we prove Theorem \ref{Theorem 2} by  using the method of blow-up analysis.

\subsection{Elementary properties of $u_\epsilon$}

In view of the equation (\ref{Euler-Lagrange}), it is important to know whether $\lambda_\epsilon$ has a positive lower bound  or not.
For this purpose, we have the following:

\begin{lemma}\label{Lemma 4}
Let $\lambda_\epsilon$ be as in (\ref{Euler-Lagrange}). Then
there holds $\liminf_{\epsilon\ra 0}\lambda_\epsilon>0$.
\end{lemma}

 {\it Proof.} For any $u\in W^{1,2}(\mathbb{R}^2)$ with $\|u\|_{1,\tau}\leq 1$, we calculate by employing Theorem \ref{prop3},
 $$\int_{\mathbb{R}^2}
\f{e^{4\pi(1-\beta)u^2}-1}{|x|^{2\beta}}dx=\lim_{\epsilon\ra 0}\int_{\mathbb{R}^2}
\f{e^{4\pi(1-\beta-\epsilon)u^2}-1}{|x|^{2\beta}}dx\leq\lim_{\epsilon\ra 0}\int_{\mathbb{R}^2}
\f{e^{4\pi(1-\beta-\epsilon)u_\epsilon^2}-1}{|x|^{2\beta}}dx.$$
This leads to
\be\label{l-t1}\Lambda_{2,\beta,\tau}=\sup_{u\in W^{1,2}(\mathbb{R}^2),\,\|u\|_{1,\tau}\leq 1}\int_{\mathbb{R}^2}
\f{e^{4\pi(1-\beta)u^2}-1}{|x|^{2\beta}}dx\leq \lim_{\epsilon\ra 0}\int_{\mathbb{R}^2}
\f{e^{4\pi(1-\beta-\epsilon)u_\epsilon^2}-1}{|x|^{2\beta}}dx.\ee
But one can easily see that
\be\label{l-t2}\int_{\mathbb{R}^2}
\f{e^{4\pi(1-\beta-\epsilon)u_\epsilon^2}-1}{|x|^{2\beta}}dx\leq\Lambda_{2,\beta,\tau}.\ee
Moreover, using the inequality $e^t\leq 1+te^t$ for $t\geq 0$, we get
$$\lambda_\epsilon\geq \f{1}{4\pi(1-\beta-\epsilon)}\int_{\mathbb{R}^2}
\f{e^{4\pi(1-\beta-\epsilon)u_\epsilon^2}-1}{|x|^{2\beta}}dx.$$
This together with (\ref{l-t1}) and (\ref{l-t2}) leads to
$$\liminf_{\epsilon\ra 0}\lambda_\epsilon\geq \lim_{\epsilon\ra 0}\f{1}{4\pi(1-\beta-\epsilon)}\int_{\mathbb{R}^2}
\f{e^{4\pi(1-\beta-\epsilon)u_\epsilon^2}-1}{|x|^{2\beta}}dx=\f{\Lambda_{2,\beta,\tau}}{4\pi(1-\beta)}>0.$$
This ends the proof of the lemma. $\hfill\Box$\\

Denote $c_\epsilon=u_\epsilon(0)=\max_{\mathbb{R}^2}u_\epsilon$. If $c_\epsilon$ is bounded, then applying elliptic estimates
to (\ref{Euler-Lagrange}), we can find some $u^\ast\in W^{1,2}(\mathbb{R}^2)$ such that $u_\epsilon\ra u^\ast$ in
$C^1_{\rm loc}(\mathbb{R}^2\setminus\{0\})\cap C^0_{\rm loc}(\mathbb{R}^2)$. Clearly $u^\ast$ is the desired
extremal function satisfying
$$\int_{\mathbb{R}^2}\f{e^{4\pi(1-\beta){u^\ast}^2}-1}{|x|^{2\beta}}dx=
\sup_{u\in W^{1,2}(\mathbb{R}^2),\,\|u\|_{1,\tau}\leq 1}\int_{\mathbb{R}^2}
\f{e^{4\pi(1-\beta)u^2}-1}{|x|^{2\beta}}dx.$$
Hence the proof of Theorem \ref{Theorem 2} terminates. In the following, we assume $c_\epsilon\ra+\infty$.
Since $u_\epsilon$ is bounded in $W^{1,2}(\mathbb{R}^2)$, we can assume without loss of generality,
$u_\epsilon$ converges to $u_0$ weakly in $W^{1,2}(\mathbb{R}^2)$, strongly in $L^q_{\rm loc}(\mathbb{R}^2)$
for any $q>0$, and a.e. in $\mathbb{R}^2$. Then we have the following:

\begin{lemma}\label{Lemma 5}
$u_0\equiv 0$ and $|\nabla u_\epsilon|^2dx\rightharpoonup \delta_0$ weakly in the sense of measure,
where $\delta_0$ denotes the Dirac measure centered at
$0\in\mathbb{R}^2$. Moreover, $u_\epsilon\ra 0$ strongly in $L^p(\mathbb{R}^2)$ for all $p\geq 2$.
\end{lemma}

{\it Proof.} For any $a\geq 0$, $b\geq 0$ and $p\geq 1$, there holds (\cite{YangJFA}, Lemma 2.1)
\be\label{ep}(e^a-1)^p\leq e^{pa}-1.\ee
Using $e^{a+b}-1=(e^a-1)(e^b-1)+(e^a-1)+(e^b-1)$, the H\"older inequality and (\ref{ep}),
we estimate
\bea\nonumber\int_{\mathbb{R}^2}\f{e^{\beta_\epsilon pu_\epsilon^2}-1}{|x|^{2\beta p}}dx&\leq&
\int_{\mathbb{R}^2}\f{e^{\beta_\epsilon p\le((1+\nu)(u_\epsilon-u_0)^2+(1+\nu^{-1})u_0^2\ri)}-1}{|x|^{2\beta p}}dx\\\nonumber
&=&\int_{\mathbb{R}^2}\f{(e^{\beta_\epsilon p(1+\nu)(u_\epsilon-u_0)^2}-1)(e^{\beta_\epsilon p(1+\nu^{-1})u_0^2}-1)}
{|x|^{2\beta p}}dx\\\nonumber
&&+\int_{\mathbb{R}^2}\f{e^{\beta_\epsilon p(1+\nu)(u_\epsilon-u_0)^2}-1}
{|x|^{2\beta p}}dx+\int_{\mathbb{R}^2}\f{e^{\beta_\epsilon p(1+\nu^{-1})u_0^2}-1}
{|x|^{2\beta p}}dx\\\nonumber
&\leq&\le(\int_{\mathbb{R}^2}\f{e^{\beta_\epsilon pp_1(1+\nu)(u_\epsilon-u_0)^2}-1}
{|x|^{2\beta p}}dx\ri)^{{1}/{p_1}}\le(\int_{\mathbb{R}^2}\f{e^{\beta_\epsilon p p_2(1+\nu^{-1})u_0^2}-1}
{|x|^{2\beta p}}dx\ri)^{{1}/{p_2}}\\ \label{Lions}
&&+\int_{\mathbb{R}^2}\f{e^{\beta_\epsilon p(1+\nu)(u_\epsilon-u_0)^2}-1}
{|x|^{2\beta p}}dx+\int_{\mathbb{R}^2}\f{e^{\beta_\epsilon p(1+\nu^{-1})u_0^2}-1}
{|x|^{2\beta p}}dx,
\eea
where $\beta_\epsilon=4\pi(1-\beta-\epsilon)$, $p>1$, $\nu>0$, $p_1>1$ and $1/p_1+1/p_2=1$.
We first prove that $u_0\equiv 0$. Suppose not. Since
\bna\|u_\epsilon-u_0\|_{1,\tau}^2=\|u_\epsilon\|_{1,\tau}^2+\|u_0\|_{1,\tau}^2-2\int_{\mathbb{R}^2}
(\nabla u_\epsilon\nabla u_0+\tau u_\epsilon u_0)dx
=1-\|u_0\|_{1,\tau}^2+o_\epsilon(1),
\ena
one can choose $p$, $p_1$ sufficiently close to $1$ and $\nu$ sufficiently close to $0$ such that
$$\f{\beta_\epsilon pp_1(1+\nu)\|u_\epsilon-u_0\|_{1,\tau}^2}{4\pi}+\f{2\beta p}{2}<1.$$
In view of the singular Trudinger-Moser inequality (\ref{Singular-RN}), we conclude that all integrals
on the right hand side of (\ref{Lions}) are bounded. Therefore
$$\int_{\mathbb{R}^2}\f{e^{4\pi(1-\beta-\epsilon)pu_\epsilon^2}-1}{|x|^{2\beta p}}dx\leq C$$
for some constant $C$ depending only on $\beta$ and $p$.
It follows that
$e^{4\pi(1-\beta-\epsilon){u_\epsilon}^2}/|x|^{2\beta}$ is bounded in $L^p(B_1)$.
This together with Lemma \ref{Lemma 4}
and $u_\epsilon$ is bounded in $L^q(B_1)$ for all $q>0$ implies that $\Delta u_\epsilon$ is bounded in $L^{p^\prime}(B_1)$ for some $p^\prime>1$.
Applying elliptic estimate to (\ref{Euler-Lagrange}), we conclude that $u_\epsilon$ is uniformly bounded in $B_{1/2}$.
This contradicts $c_\epsilon\ra+\infty$. Therefore
$u_0\equiv 0$.

We next prove that $|\nabla u_\epsilon|^2dx\rightharpoonup \delta_0$. For otherwise, we can choose sufficiently small $\bar{r}>0$ such that
$$\int_{B_{\bar{r}}}(|\nabla u_\epsilon|^2+\tau u_\epsilon^2)dx\leq \eta<1$$
for sufficiently small $\epsilon>0$. Hence $\Delta u_\epsilon$ is bounded in $L^{p^{\prime\prime}}(B_{\bar{r}})$
for some $p^{\prime\prime}>1$, and thus
elliptic estimate leads to  $u_\epsilon$ is uniformly bounded in $B_{\bar{r}/2}$ contradicting $c_\epsilon\ra+\infty$.

For the last assertion, noting that $\|u_\epsilon\|_{1,\tau}=1$ and $|\nabla u_\epsilon|^2dx\rightharpoonup \delta_0$,
we obtain $\|u_\epsilon\|_{L^2(\mathbb{R}^2)}=o_\epsilon(1)$. Taking $M>0$ such that if $|x|>M$, then $u_\epsilon<1$, one has for
any $p>2$,
\bna
\int_{\mathbb{R}^2}u_\epsilon^pdx&=&\int_{|x|>M}u_\epsilon^pdx+\int_{|x|\leq M}u_\epsilon^pdx\\
&\leq& \int_{|x|>M}u_\epsilon^2dx+o_\epsilon(1)\\
&\leq& \int_{\mathbb{R}^2}u_\epsilon^2dx+o_\epsilon(1)=o_\epsilon(1).\ena
Here we have used the fact that $u_\epsilon\ra 0$ in $L^q_{\rm loc}(\mathbb{R}^2)$ for any $q>0$. This completes the proof of the lemma.
 $\hfill\Box$

 \subsection{Blow-up analysis}

Set $r_\epsilon=\sqrt{\lambda_\epsilon}c_\epsilon^{-1}e^{-2\pi(1-\beta-\epsilon)c_\epsilon^2}$,
$\psi_\epsilon(x)=c_\epsilon^{-1}
u_\epsilon(r_\epsilon^{{1}/{(1-\beta)}}x)$ and $\varphi_\epsilon(x)=c_\epsilon(u_\epsilon(r_\epsilon^{{1}/{(1-\beta)}}x)-c_\epsilon)$.
Then we have the following:

\begin{lemma}\label{Lemma 6} $(i)$ For any $\gamma<2\pi(1-\beta)$, there holds
$r_\epsilon e^{\gamma c_\epsilon^2}\ra 0$ as $\epsilon\ra 0$; $(ii)$ $\psi_\epsilon\ra 1$ in $C^1_{\rm loc}(\mathbb{R}^2\setminus\{0\})\cap C^0_{\rm loc}(\mathbb{R}^2)$; $(iii)$
$\varphi_\epsilon\ra \varphi$ in
$C^1_{\rm loc}(\mathbb{R}^2\setminus\{0\})\cap C^0_{\rm loc}(\mathbb{R}^2)$, where
$\varphi(x)=-\f{1}{4\pi(1-\beta)}\log(1+\f{\pi}{1-\beta}|x|^{2(1-\beta)})$ and
$\int_{\mathbb{R}^2}|x|^{-2\beta}e^{8\pi(1-\beta)\varphi}dx=1$.
\end{lemma}

{\it Proof.} $(i)$  By definition of $r_\epsilon$, we have
\bea\nonumber
r_\epsilon^2e^{2\gamma c_\epsilon^2}&=&c_\epsilon^{-2}e^{-4\pi(1-\beta-\epsilon-\f{\gamma}{2\pi})c_\epsilon^2}
\int_{\mathbb{R}^2}\f{u_\epsilon^2e^{4\pi(1-\beta-\epsilon)u_\epsilon^2}}{|x|^{2\beta}}dx\\\nonumber
&\leq& c_\epsilon^{-2}\int_{\mathbb{R}^2}\f{u_\epsilon^2e^{2\gamma u_\epsilon^2}}{|x|^{2\beta}}dx
\\\label{tds}
&\leq& c_\epsilon^{-2}\int_{\mathbb{R}^2}\f{u_\epsilon^2(e^{2\gamma u_\epsilon^2}-1)}{|x|^{2\beta}}dx
+c_\epsilon^{-2}\int_{\mathbb{R}^2}\f{u_\epsilon^2}{|x|^{2\beta}}dx.
\eea
By Lemma \ref{Lemma 5}, we know that $\|u_\epsilon\|_{L^p(\mathbb{R}^2)}=o_\epsilon(1)$ for any $p\geq 2$.
As an easy consequence, there holds for any $p\geq 2$
\be\label{lim-0}\lim_{\epsilon\ra 0}\int_{\mathbb{R}^2}\f{u_\epsilon^{p}}{|x|^{2\beta}}dx=0.\ee
Noting that $\gamma<2\pi(1-\beta)$, we can choose $p_1>1$ such that $\gamma p_1<2\pi(1-\beta)$. In view of
(\ref{Singular-RN}), (\ref{ep}) and (\ref{lim-0}), we have by the H\"older inequality,
\be\label{t-2}\int_{\mathbb{R}^2}\f{u_\epsilon^2(e^{2\gamma u_\epsilon^2}-1)}{|x|^{2\beta}}dx\leq
\le(\int_{\mathbb{R}^2}\f{e^{2\gamma p_1 u_\epsilon^2}-1}{|x|^{2\beta}}dx\ri)^{1/p_1}
\le(\int_{\mathbb{R}^2}\f{u_\epsilon^{2p_2}}{|x|^{2\beta}}dx\ri)^{1/p_2}=o_\epsilon(1),\ee
where $1/p_2+1/p_1=1$. Inserting (\ref{lim-0}) and (\ref{t-2}) into (\ref{tds}), we obtain
$r_\epsilon e^{\gamma c^2_\epsilon}\ra 0$ as $\epsilon\ra 0$.

$(ii)$ It can be easily checked that $\psi_\epsilon$ satisfies the equation
\be\label{psi-e}-\Delta\psi_\epsilon(x)=-\tau r_\epsilon^{{2}/{(1-\beta)}}\psi_\epsilon(x)+c_\epsilon^{-2}
|x|^{-2\beta}\psi_\epsilon(x) e^{4\pi(1-\beta-\epsilon)(u_\epsilon^2(r_\epsilon^{1/(1-\beta)}x)-c_\epsilon^2)}.\ee
Since $|\psi_\epsilon|\leq 1$, $u_\epsilon^2\leq c_\epsilon^2$ and $r_\epsilon\ra 0$ as $\epsilon\ra 0$,
we have by applying elliptic estimates to (\ref{psi-e}), $\psi_\epsilon\ra \psi$ in $C^1_{\rm loc}(\mathbb{R}^2\setminus\{0\})
\cap C^0_{\rm loc}(\mathbb{R}^2)$, where $\psi$ is a bounded harmonic function on $\mathbb{R}^2$. Then the Liouville theorem
leads to $\psi\equiv\psi(0)=1$.

$(iii)$ A straightforward calculation shows
\be\label{varp}-\Delta\varphi_\epsilon(x)=-\tau c_\epsilon^2r_\epsilon^{{2}/{(1-\beta)}}\psi_\epsilon(x)+
|x|^{-2\beta}\psi_\epsilon(x)e^{4\pi(1-\beta-\epsilon)(1+\psi_\epsilon(x))\varphi_\epsilon(x)}.\ee
Note that $\varphi_\epsilon(x)\leq 0=\max_{\mathbb{R}^2}\varphi_\epsilon$. Applying elliptic estimates to (\ref{varp}),
we conclude that $\varphi_\epsilon\ra \varphi$ in $C^1_{\rm loc}(\mathbb{R}^2\setminus\{0\})
\cap C^0_{\rm loc}(\mathbb{R}^2)$, where $\varphi$ is a distributional solution to
\be\label{vpi}\le\{\begin{array}{lll}-\Delta \varphi(x)=|x|^{-2\beta}e^{8\pi(1-\beta)\varphi(x)}\quad{\rm in}\quad \mathbb{R}^2,\\
[1.5ex] \varphi(0)=0.\end{array}\ri.\ee
Since $u_\epsilon$ is decreasing symmetric and $u_\epsilon(0)=\max_{\mathbb{R}^2}u_\epsilon=c_\epsilon$, $\varphi$ must be
decreasing symmetric and $\varphi(0)=\max_{\mathbb{R}^2}\varphi$. If we set $\bar\varphi(r)=\varphi(x)$ for any $x\in \mathbb{R}^2$
and $r=|x|$, then (\ref{vpi}) reduces to
\be\label{vpi-rad}\le\{\begin{array}{lll}-(r\bar{\varphi}^\prime)^\prime=r^{1-2\beta}e^{8\pi(1-\beta)\bar\varphi},\\
[1.5ex] \bar\varphi(0)=0.\end{array}\ri.\ee
Clearly, this equation has a special solution
$$\bar\varphi(r)=-\f{1}{4\pi(1-\beta)}\log(1+\f{\pi}{1-\beta}r^{2(1-\beta)}).$$
By the standard uniqueness result of the ordinary differential equation (\ref{vpi-rad}), we have
$$\varphi(x)=-\f{1}{4\pi(1-\beta)}\log(1+\f{\pi}{1-\beta}|x|^{2(1-\beta)}),\quad x\in\mathbb{R}^2.$$
It follows that
\be\label{energy}\int_{\mathbb{R}^2}|x|^{-2\beta}e^{8\pi(1-\beta)\varphi}dx=\int_0^{+\infty}\f{2\pi r^{1-2\beta}}{(1+\f{\pi}
{1-\beta}r^{2(1-\beta)})^2}dr=\int_0^{+\infty}\f{1}{(1+t)^2}dt=1.\ee
This completes the proof of the lemma. $\hfill\Box$\\

Lemma \ref{Lemma 6} gives convergence behavior of $u_\epsilon$ near $0$. To reveal the convergence behavior of $u_\epsilon$ away
from $0$, following \cite{Li-Ruf},  we define $u_{\epsilon,\gamma}=\min\{u_\epsilon,\gamma c_\epsilon\}$ for any $0<\gamma<1$.
Then we have the following:

\begin{lemma}\label{Lemma 7}
For any $0<\gamma<1$, there holds $\lim_{\epsilon\ra 0}\int_{\mathbb{R}^2}|\nabla u_{\epsilon,\gamma}|^2dx=\gamma$.
\end{lemma}

{\it Proof.} Testing the equation (\ref{Euler-Lagrange})
by $u_{\epsilon,\gamma}$, we have for any fixed $R>0$,
\bna
\int_{\mathbb{R}^2}|\nabla u_{\epsilon,\gamma}|^2dx&=&-\tau\int_{\mathbb{R}^2}u_\epsilon u_{\epsilon,\gamma}dx+
\f{1}{\lambda_\epsilon}\int_{\mathbb{R}^2}u_\epsilon u_{\epsilon,\gamma}\f{e^{4\pi(1-\beta-\epsilon)u_\epsilon^2}}
{|x|^{2\beta}}dx\\
&\geq&\f{1}{\lambda_\epsilon}\int_{B_{Rr_\epsilon^{1/(1-\beta)}}}\gamma c_\epsilon u_\epsilon
\f{e^{4\pi(1-\beta-\epsilon)u_\epsilon^2}}
{|x|^{2\beta}}dx+o_\epsilon(1)\\
&=&(1+o_\epsilon(1))\gamma\int_{B_R}\f{e^{8\pi(1-\beta)\varphi(x)}}{|x|^{2\beta}}dx+o_\epsilon(1).
\ena
Hence
$$\liminf_{\epsilon\ra 0}\int_{\mathbb{R}^2}|\nabla u_{\epsilon,\gamma}|^2dx\geq \gamma\int_{B_R}
\f{e^{8\pi(1-\beta)\varphi(x)}}{|x|^{2\beta}}dx.$$
In view of (\ref{energy}), passing to the limit $R\ra+\infty$, we obtain
\be\label{geq}\liminf_{\epsilon\ra 0}\int_{\mathbb{R}^2}|\nabla u_{\epsilon,\gamma}|^2dx\geq \gamma.\ee
Testing the equation (\ref{Euler-Lagrange}) by $(u_\epsilon-\gamma c_\epsilon)^+$, we obtain for any fixed $R>0$,
\bna
\int_{\mathbb{R}^2}|\nabla(u_\epsilon-\gamma c_\epsilon)^+|^2dx&=&-\tau\int_{\mathbb{R}^2}u_\epsilon(u_\epsilon-\gamma c_\epsilon)^+dx
+\int_{\mathbb{R}^2}(u_\epsilon-\gamma c_\epsilon)^+ u_\epsilon\f{e^{4\pi(1-\beta-\epsilon)u_\epsilon^2}}
{\lambda_\epsilon|x|^{2\beta}}dx\\
&\geq& \f{1}{\lambda_\epsilon}\int_{B_{Rr_\epsilon^{1/(1-\beta)}}} u_\epsilon (u_\epsilon-\gamma c_\epsilon)^+
\f{e^{4\pi(1-\beta-\epsilon)u_\epsilon^2}}
{|x|^{2\beta}}dx+o_\epsilon(1)\\
&=&(1+o_\epsilon(1))(1-\gamma)\int_{B_R}\f{e^{8\pi(1-\beta)\varphi(x)}}{|x|^{2\beta}}dx+o_\epsilon(1).
\ena
Similarly as above, we have
\be\label{leq}\liminf_{\epsilon\ra 0}\int_{\mathbb{R}^2}|\nabla (u_\epsilon-\gamma c_\epsilon)^+|^2dx\geq 1-\gamma.\ee
Note that
\be\label{1}\int_{\mathbb{R}^2}|\nabla u_{\epsilon,\gamma}|^2dx+
\int_{\mathbb{R}^2}|\nabla (u_\epsilon-\gamma c_\epsilon)^+|^2dx=\int_{\mathbb{R}^2}|\nabla u_{\epsilon}|^2dx=1+o_\epsilon(1).\ee
Combining (\ref{geq}), (\ref{leq}) and (\ref{1}), we conclude the lemma. $\hfill\Box$

\begin{lemma}\label{Lemma 8}
There holds
\be\label{lit}\lim_{\epsilon\ra 0}\int_{\mathbb{R}^2}\f{e^{4\pi(1-\beta-\epsilon)u_\epsilon^2}-1}{|x|^{2\beta}}dx=\lim_{\epsilon\ra 0}
\f{\lambda_\epsilon}{c_\epsilon^2}.\ee
\end{lemma}

{\it Proof.} Let $0<\gamma<1$ be fixed. Using the inequality $e^t-1\leq te^t$ $(t\geq 0)$ and the definition of
$u_{\epsilon,\gamma}$, we obtain
\bea\nonumber
\int_{u_\epsilon\leq \gamma c_\epsilon}\f{e^{4\pi(1-\beta-\epsilon)u_\epsilon^2}-1}{|x|^{2\beta}}dx&\leq&\int_{\mathbb{R}^2}
\f{e^{4\pi(1-\beta-\epsilon)u_{\epsilon,\gamma}^2}-1}{|x|^{2\beta}}dx\\\nonumber
&\leq&4\pi(1-\beta)\int_{\mathbb{R}^2}
\f{u_{\epsilon,\gamma}^2e^{4\pi(1-\beta-\epsilon)u_{\epsilon,\gamma}^2}}{|x|^{2\beta}}dx\\\label{gam1}
&=& 4\pi(1-\beta)\le\{\int_{\mathbb{R}^2}u_{\epsilon,\gamma}^2
\f{e^{4\pi(1-\beta-\epsilon)u_{\epsilon,\gamma}^2}-1}{|x|^{2\beta}}dx+\int_{\mathbb{R}^2}
\f{u_{\epsilon,\gamma}^2}{|x|^{2\beta}}dx\ri\}.
\eea
It follows from (\ref{lim-0}) that
\be\label{gam1-1}\int_{\mathbb{R}^2}
\f{u_{\epsilon,\gamma}^2}{|x|^{2\beta}}dx\leq \int_{\mathbb{R}^2}
\f{u_{\epsilon}^2}{|x|^{2\beta}}dx=o_\epsilon(1).\ee
Moreover, combining Lemma \ref{Lemma 5} and Lemma \ref{Lemma 7}, we have
$\lim_{\epsilon\ra 0}\|u_{\epsilon,\gamma}\|_{1,\tau}^2=\gamma<1$. Let
$1<p<1/\gamma$ be fixed and $1/p+1/p^\prime=1$. Using the H\"older inequality and
the singular Trudinger-Moser inequality (\ref{Singular-RN}), we have
\bea\nonumber
\int_{\mathbb{R}^2}u_{\epsilon,\gamma}^2
\f{e^{4\pi(1-\beta-\epsilon)u_{\epsilon,\gamma}^2}-1}{|x|^{2\beta}}dx&\leq&
\le(\int_{\mathbb{R}^2}
\f{e^{4\pi(1-\beta-\epsilon)pu_{\epsilon,\gamma}^2}-1}{|x|^{2\beta}}dx\ri)^{1/p}
\le(\int_{\mathbb{R}^2}
\f{u_{\epsilon,\gamma}^{2p^\prime}}{|x|^{2\beta}}dx\ri)^{1/p^\prime}\\\label{gam1-2}
&\leq& C\le(\int_{\mathbb{R}^2}
\f{u_{\epsilon,\gamma}^{2p^\prime}}{|x|^{2\beta}}dx\ri)^{1/p^\prime}
\eea
for some constant $C$ depending only on $\beta$, $p$ and $\gamma$. Inserting (\ref{gam1-1}) and (\ref{gam1-2}) into
(\ref{gam1}), one has
\be\label{gam1-00}
\int_{u_\epsilon\leq \gamma c_\epsilon}\f{e^{4\pi(1-\beta-\epsilon)u_\epsilon^2}-1}{|x|^{2\beta}}dx=o_\epsilon(1).
\ee
Moreover, we estimate
\bea\nonumber
\int_{u_\epsilon> \gamma c_\epsilon}\f{e^{4\pi(1-\beta-\epsilon)u_\epsilon^2}-1}{|x|^{2\beta}}dx&=&
\int_{u_\epsilon> \gamma c_\epsilon}\f{e^{4\pi(1-\beta-\epsilon)u_\epsilon^2}}{|x|^{2\beta}}dx+o_\epsilon(1)\\
\nonumber&\leq&\f{1}{\gamma^2}\int_{u_\epsilon> \gamma c_\epsilon}\f{u_\epsilon^2}{c_\epsilon^2}
\f{e^{4\pi(1-\beta-\epsilon)u_\epsilon^2}}{|x|^{2\beta}}dx+o_\epsilon(1)\\\label{gam2}
&\leq& \f{1}{\gamma^2}\f{\lambda_\epsilon}{c_\epsilon^2}+o_\epsilon(1).
\eea
Combining (\ref{gam1-00}) and (\ref{gam2}), we have
$$\lim_{\epsilon\ra 0}\int_{\mathbb{R}^2}\f{e^{4\pi(1-\beta-\epsilon)u_\epsilon^2}-1}{|x|^{2\beta}}dx\leq
\f{1}{\gamma^2}\liminf_{\epsilon\ra 0}\f{\lambda_\epsilon}{c_\epsilon^2}.$$
Letting $\gamma\ra 1$, we conclude
\be\label{ll}
\lim_{\epsilon\ra 0}\int_{\mathbb{R}^2}\f{e^{4\pi(1-\beta-\epsilon)u_\epsilon^2}-1}{|x|^{2\beta}}dx\leq
\liminf_{\epsilon\ra 0}\f{\lambda_\epsilon}{c_\epsilon^2}.
\ee
On the other hand,
\bna
\f{\lambda_\epsilon}{c_\epsilon^2}&=&\int_{\mathbb{R}^2}\f{u_\epsilon^2}{c_\epsilon^2}
\f{e^{4\pi(1-\beta-\epsilon)u_\epsilon^2}}{|x|^{2\beta}}dx\\
&=&\int_{\mathbb{R}^2}\f{u_\epsilon^2}{c_\epsilon^2}
\f{e^{4\pi(1-\beta-\epsilon)u_\epsilon^2}-1}{|x|^{2\beta}}dx+\f{1}{c_\epsilon^2}\int_{\mathbb{R}^2}
\f{u_\epsilon^2}{|x|^{2\beta}}dx\\
&\leq&\int_{\mathbb{R}^2}
\f{e^{4\pi(1-\beta-\epsilon)u_\epsilon^2}-1}{|x|^{2\beta}}dx+o_\epsilon(1).
\ena
Thus
\be\label{gam3}\limsup_{\epsilon\ra 0}\f{\lambda_\epsilon}{c_\epsilon^2}\leq
\lim_{\epsilon\ra 0}\int_{\mathbb{R}^2}
\f{e^{4\pi(1-\beta-\epsilon)u_\epsilon^2}-1}{|x|^{2\beta}}dx.\ee
Combining (\ref{ll}) and (\ref{gam3}), we get the desired result. $\hfill\Box$

\begin{corollary}\label{Cor}
If $\theta<2$, then $\lambda_\epsilon/c_\epsilon^\theta\ra +\infty$ as $\epsilon\ra 0$.
\end{corollary}

{\it Proof.} An obvious consequence of Lemma \ref{Lemma 8}. $\hfill\Box$

\begin{lemma}\label{Lemma 9}
For any function $\phi\in C_0^0(\mathbb{R}^2)$, there holds
$$\lim_{\epsilon\ra 0}\int_{\mathbb{R}^2}\f{c_\epsilon u_\epsilon}{\lambda_\epsilon}
\f{e^{4\pi(1-\beta-\epsilon)u_\epsilon^2}}{|x|^{2\beta}}\phi dx=\phi(0).$$
\end{lemma}

{\it Proof.} Let $\phi\in C_0^0(\mathbb{R}^2)$ be fixed.
Write for simplicity $h_\epsilon=\lambda_\epsilon^{-1}|x|^{-2\beta}c_\epsilon u_\epsilon e^{4\pi(1-\beta-\epsilon)u_\epsilon^2}$.
Given $0<\gamma<1$. Firstly we calculate
\bna
\int_{u_\epsilon\leq\gamma c_\epsilon}h_\epsilon \phi dx=\f{c_\epsilon}{\lambda_\epsilon}\int_{u_\epsilon\leq\gamma c_\epsilon}
u_\epsilon \phi \f{e^{4\pi(1-\beta-\epsilon)u_\epsilon^2}-1}{|x|^{2\beta}}dx+\f{c_\epsilon}{\lambda_\epsilon}\int_{u_\epsilon\leq\gamma c_\epsilon}
 \f{u_\epsilon \phi}{|x|^{2\beta}}dx.
\ena
 In view of an obvious analog of (\ref{gam1-2}), there holds
\bna
\le|\int_{u_\epsilon\leq\gamma c_\epsilon}
u_\epsilon \phi \f{e^{4\pi(1-\beta-\epsilon)u_\epsilon^2}-1}{|x|^{2\beta}}dx\ri|&\leq&
\le(\sup_{\mathbb{R}^2}|\phi|\ri)\int_{\mathbb{R}^2}u_{\epsilon,\gamma}\f{e^{4\pi(1-\beta-\epsilon)
u_{\epsilon,\gamma}^2}-1}{|x|^{2\beta}}dx=o_\epsilon(1).
\ena
Note that $u_\epsilon\ra 0$ in $L^q_{\rm loc}(\mathbb{R}^2)$ for any $q>0$. We derive
$$\le|\int_{u_\epsilon\leq\gamma c_\epsilon}
 \f{u_\epsilon \phi}{|x|^{2\beta}}dx\ri|\leq\le(\sup_{\mathbb{R}^2}|\phi|\ri)\int_{{\rm supp}\,\phi}
 \f{u_\epsilon}{|x|^{2\beta}}dx=o_\epsilon(1).$$
 By Corollary \ref{Cor}, $c_\epsilon/\lambda_\epsilon=o_\epsilon(1)$. Therefore
 \be\label{g1}\int_{u_\epsilon\leq\gamma c_\epsilon}h_\epsilon \phi dx=o_\epsilon(1).\ee
 It follows from Lemma \ref{Lemma 6} that $B_{Rr_\epsilon^{1/(1-\beta)}}\subset\{u_\epsilon>\gamma c_\epsilon\}$
 for sufficiently small $\epsilon>0$, that
 \bna
 \int_{B_{Rr_\epsilon^{1/(1-\beta)}}}h_\epsilon\phi dx&=&\phi(0)(1+o_\epsilon(1))\le(\int_{B_R}
 \f{e^{8\pi(1-\beta)\varphi}}{|x|^{2\beta}}dx+o_\epsilon(1)\ri)\\
 &=&\phi(0)(1+o_\epsilon(1)+o_R(1)),
 \ena
 and that
 \bna
 \le|\int_{\{u_\epsilon>\gamma c_\epsilon\}\setminus B_{Rr_\epsilon^{1/(1-\beta)}}}h_\epsilon\phi dx\ri|
 &\leq&\f{1}{\gamma}\le(\sup_{\mathbb{R}^2}|\phi|\ri)\int_{\{u_\epsilon>\gamma c_\epsilon\}\setminus B_{Rr_\epsilon^{1/(1-\beta)}}}
 \f{u_\epsilon^2}{\lambda_\epsilon}\f{e^{4\pi(1-\beta-\epsilon)u_\epsilon^2}}{|x|^{2\beta}}dx\\
 &\leq&\f{1}{\gamma}\le(\sup_{\mathbb{R}^2}|\phi|\ri)\le(1-\int_{B_R}\f{e^{8\pi(1-\beta)\varphi}}{|x|^{2\beta}}dx+o_\epsilon(1)\ri)\\
 &=& o_\epsilon(1)+o_R(1).
 \ena
 It then follows that
 \be\label{g2}\lim_{\epsilon\ra 0}\int_{u_\epsilon>\gamma c_\epsilon}h_\epsilon \phi dx=\phi(0).\ee
 Combining (\ref{g1}) and (\ref{g2}), we complete the proof of the lemma. $\hfill\Box$

 \begin{lemma}\label{Green}
 $c_\epsilon u_\epsilon \ra G$ in $C^1_{\rm loc}(\mathbb{R}^2\setminus\{0\})$ and weakly in $W^{1,q}(\mathbb{R}^2)$
 for all $1<q<2$, where $G$ is a distributional solution to
 \be\label{gr}-\Delta G+\tau G=\delta_0\quad{\rm in}\quad\mathbb{R}^2.\ee
 Moreover, $G\in W^{1,2}(\mathbb{R}^2\setminus B_r)$ for any $r>0$ and $G$ takes the form
\be\label{gh}G(x)=-\f{1}{2\pi}\log|x|+A_0+w(x),\ee
where $A_0$ is a constant, $w\in C^1(\mathbb{R}^2)$ and $w(0)=0$.
 \end{lemma}

 {\it Proof.} Multiplying both sides of the equation (\ref{Euler-Lagrange}) by $c_\epsilon$, one has
 \be\label{g-e}-\Delta (c_\epsilon u_\epsilon)+\tau (c_\epsilon u_\epsilon)=\f{c_\epsilon u_\epsilon
e^{4\pi(1-\beta-\epsilon) u_\epsilon^2}}{\lambda_\epsilon|x|^{2\beta}}\quad{\rm in}\quad \mathbb{R}^2.\ee
In view of Lemma \ref{Lemma 9}, $h_\epsilon=\lambda_\epsilon^{-1}|x|^{-2\beta}c_\epsilon u_\epsilon e^{4\pi(1-\beta-\epsilon) u_\epsilon^2}$
is bounded in $L^1_{\rm loc}(\mathbb{R}^2)$. Using an argument of Li-Ruf (\cite{Li-Ruf}, Proposition 3.7), which is adapted from
that of Struwe (\cite{Str-1}, Theorem 2.2), one concludes that $c_\epsilon u_\epsilon$ is bounded in $W^{1,q}_{\rm loc}(\mathbb{R}^2)$ for all
$1<q<2$. Hence $c_\epsilon u_\epsilon\rightharpoonup G$ weakly in $W^{1,q}_{\rm loc}(\mathbb{R}^2)$ for any $1<q<2$ and $G$ is a distributional
solution to (\ref{gr}). Since $\Delta(G(x)+\f{1}{2\pi}\log |x|)\in L^p_{\rm loc}(\mathbb{R}^2)$ for any $p>2$,
(\ref{gh}) follows from elliptic estimates immediately.
  Applying
elliptic estimates to the equation (\ref{g-e}), we obtain $c_\epsilon u_\epsilon\ra G$ in $C^1_{\rm loc}(\mathbb{R}^2\setminus\{0\})$.
Note that $c_\epsilon u_\epsilon\in W^{1,2}(\mathbb{R}^2)$. Multiplying both sides of (\ref{g-e}) by $c_\epsilon u_\epsilon$ and integrating
by parts on the domain $\mathbb{R}^2\setminus B_r$ for some $r>0$, we get
\bna
\int_{\mathbb{R}^2\setminus B_r}(|\nabla (c_\epsilon u_\epsilon)|^2+\tau (c_\epsilon u_\epsilon)^2)dx&=&
-\int_{\p B_r}c_\epsilon u_\epsilon\f{\p(c_\epsilon u_\epsilon)}{\p \nu}d\sigma+\int_{\mathbb{R}^2\setminus B_r}h_\epsilon c_\epsilon u_\epsilon dx\\
&\leq&-\int_{\p B_r}c_\epsilon u_\epsilon\f{\p(c_\epsilon u_\epsilon)}{\p \nu}d\sigma+\f{c_\epsilon^2e^{4\pi u_\epsilon^2(r)}}{\lambda_\epsilon}
\int_{\mathbb{R}^2\setminus B_r}\f{u_\epsilon^2}{|x|^{2\beta}}dx\\
&\leq& C_r
\ena
for some constant $C_r$ depending only on $r$, since $c_\epsilon u_\epsilon\ra G$ in $C^1_{\rm loc}(\mathbb{R}^2\setminus\{0\})$. This also leads to
$$\int_{r\leq |x|\leq R}(|\nabla G|^2+\tau G^2)dx\leq C_r, \quad\forall R>r.$$
Passing to the limit $R\ra\infty$, we have
$$\int_{\mathbb{R}^2\setminus B_r}(|\nabla G|^2+\tau G^2)dx\leq C_r.$$
This gives the desired result. $\hfill\Box$

\subsection{Upper bound estimate}

We need a singular version of Carleson-Chang's upper bound estimate, namely Lemma \ref{C-C-Lemma}.

\begin{lemma}{\label{Lemma 10}}
Let $w_\epsilon\in W_0^{1,2}(B_r)$ satisfies $\int_{B_r}|\nabla w_\epsilon|^2dx\leq 1$, $w_\epsilon\rightharpoonup 0$
weakly in $W_0^{1,2}(B_r)$, and $w_\epsilon$ is radially symmetric. Then
\be\label{r}\limsup_{\epsilon\ra 0}\int_{B_r}\f{e^{4\pi(1-\beta)w_\epsilon^2}-1}{|x|^{2\beta}}dx\leq \f{e\pi}{1-\beta}r^{2(1-\beta)}.\ee
\end{lemma}

{\it Proof.} We first prove (\ref{r}) for $r=1$.

Denote $w_\epsilon(|x|)=w_\epsilon(x)$. Let $v_\epsilon(x)=\sqrt{1-\beta}w_\epsilon(|x|^{{1}/{(1-\beta)}})$. Then
$$\int_{B_1}|\nabla v_\epsilon|^2dx=\int_{B_1}|\nabla w_\epsilon|^2dx.$$
Clearly we can assume up to a subsequence, $v_\epsilon\rightharpoonup v_0$ weakly in $W_0^{1,2}(B_1)$, $v_\epsilon\ra v_0$
strongly in $L^2(B_1)$, and $v_\epsilon\ra v_0$ a.e. in $B_1$. Also, we can assume $w_\epsilon\ra 0$ a.e. in $B_1$. Hence
we conclude $v_0=0$ a.e. in $B_1$. By a change of variable $t=s^{1/(1-\beta)}$, there holds
\bna
\int_{B_1}\f{e^{4\pi(1-\beta)w_\epsilon^2}-1}{|x|^{2\beta}}dx&=&\int_0^1\f{e^{4\pi(1-\beta)w_\epsilon^2(t)}-1}{t^{2\beta}}
2\pi tdt\\
&=&\f{2\pi}{1-\beta}\int_0^1 s^{{(1-2\beta)}/{(1-\beta)}}(e^{4\pi(1-\beta)w_\epsilon^2(s^{1/(1-\beta)})}-1)s^{\beta/(1-\beta)}ds\\
&=&\f{2\pi}{1-\beta}\int_0^1  s(e^{4\pi v_\epsilon^2(s)}-1)d s\\
&=&\f{1}{1-\beta}\int_{B_1}(e^{4\pi v_\epsilon^2}-1)dx.
\ena
It follows from Lemma \ref{C-C-Lemma} that
\be\label{bound}\limsup_{\epsilon\ra 0}\int_{B_1}\f{e^{4\pi(1-\beta)w_\epsilon^2}-1}{|x|^{2\beta}}dx
\leq \f{e\pi}{1-\beta}.\ee

We next prove (\ref{r}) for the case of general $r$. Set $\tilde{w}_\epsilon(x)=w_\epsilon(rx)$ for $x\in B_1$.
One can check that
$$\int_{B_1}|\nabla \tilde{w}_\epsilon|^2dx=\int_{B_r}|\nabla w_\epsilon|^2dx$$
and that
\bna
\int_{B_r}\f{e^{4\pi(1-\beta)w_\epsilon^2}-1}{|x|^{2\beta}}dx=r^{2(1-\beta)}\int_{B_1}\f{e^{4\pi(1-\beta)\tilde{w}_\epsilon^2}-1}
{|x|^{2\beta}}dx.
\ena
This together with (\ref{bound}) gives the desired result. $\hfill\Box$

By the equation (\ref{Euler-Lagrange}) and $\|u_\epsilon\|_{1,\tau}=1$, we have
\bea\nonumber
\int_{B_r}|\nabla u_\epsilon|^2dx&=&1-\int_{\mathbb{R}^2\setminus B_r}(|\nabla u_\epsilon|^2+\tau u_\epsilon^2)dx
-\tau\int_{B_r}u_\epsilon^2dx\\\label{shang}
&=&1-\int_{\mathbb{R}^2\setminus B_r}\f{u_\epsilon^2}{\lambda_\epsilon}\f{e^{4\pi(1-\beta-\epsilon)u_\epsilon^2}}{|x|^{2\beta}}dx
+\int_{\p B_r}u_\epsilon\f{\p u_\epsilon}{\p r}d\sigma-\tau\int_{B_r}u_\epsilon^2dx.
\eea
Since
\bna\int_{\mathbb{R}^2\setminus B_r}\f{u_\epsilon^2}{\lambda_\epsilon}\f{e^{4\pi(1-\beta-\epsilon)u_\epsilon^2}}{|x|^{2\beta}}dx&=&
\f{1}{c_\epsilon^2}\f{c_\epsilon^2}{\lambda_\epsilon}\int_{\mathbb{R}^2\setminus B_r}{u_\epsilon^2}
\f{e^{4\pi(1-\beta-\epsilon)u_\epsilon^2}}{|x|^{2\beta}}dx\\
&=&\f{o_\epsilon(1)}{c_\epsilon^2},
\ena
$$\int_{\p B_r}u_\epsilon\f{\p u_\epsilon}{\p r}d\sigma=\f{1}{c_\epsilon^2}\le(\int_{\p B_r}G\f{\p G}{\p r}d\sigma+o_\epsilon(1)\ri),$$
and
$$\int_{B_r}u_\epsilon^2dx=\f{1}{c_\epsilon^2}\le(\int_{B_r}G^2dx+o_\epsilon(1)\ri).$$
Inserting these equations into (\ref{shang}) and noting that $G(x)=-\f{1}{2\pi}\log|x|+A_0+w(x)$, we conclude
\be\label{r-energy}\int_{B_r}|\nabla u_\epsilon|^2dx=1-\f{1}{c_\epsilon^2}\le(\f{1}{2\pi}\log \f{1}{r}+A_0+o_\epsilon(1)+o_r(1)\ri).\ee

Denote $s_{\epsilon,r}=\sup_{\p B_r}u_\epsilon=u_\epsilon(r)$ and
          ${u}_{\epsilon,r}=(u_\epsilon-s_{\epsilon,r})^+$, the positive part of $u_\epsilon-s_{\epsilon,r}$.
          Clearly we have
          ${u}_{\epsilon,r}\in W_0^{1,2}(B_r)$.
          In view of Lemma \ref{Lemma 10},
          \be\label{r2}\limsup_{\epsilon\ra 0}\int_{B_r}\f{e^{4\pi(1-\beta){u}_{\epsilon,r}^2/\tau_{\epsilon,r}}-1}{|x|^{2\beta}}dx\leq \f{e\pi}{1-\beta}r^{2(1-\beta)},\ee
          where $\tau_{\epsilon,r}=\int_{B_r}|\nabla u_\epsilon|^2dx$. Moreover, we know from Lemma \ref{Lemma 6} that
          $u_\epsilon=c_\epsilon+o_\epsilon(1)$ on $B_{Rr_\epsilon^{1/(1-\beta)}}$. Hence, in view of (\ref{r-energy}), there holds on
          $B_{Rr_\epsilon^{1/(1-\beta)}}\subset B_r$,
          \bna
          4\pi(1-\beta-\epsilon)
          u_\epsilon^2&\leq&4\pi(1-\beta)({u}_{\epsilon,r}+s_{\epsilon,r})^2\\
          &=&
          4\pi(1-\beta){u}_{\epsilon,r}^2+8\pi(1-\beta)
          s_{\epsilon,r}{u}_{\epsilon,r}+o(1)\\
          &=&4\pi(1-\beta){u}_{\epsilon,r}^2-4(1-\beta)\log r+8\pi(1-\beta)
          A_{0}+o(1)\\
          &=&4\pi(1-\beta){u}_{\epsilon,r}^2/\tau_{\epsilon,r}-2(1-\beta)\log r+4\pi(1-\beta)
          A_{0}+o(1),
          \ena
          where $o(1)\ra 0$ as $\epsilon\ra 0$ first and next $r\ra 0$.
          Therefore
          \bea\nonumber
          \int_{B_{Rr_\epsilon^{1/(1-\beta)}}}\f{e^{4\pi(1-\beta-\epsilon)
          u_\epsilon^2}-1}{|x|^{2\beta}}dx&\leq& r^{-2(1-\beta)}e^{4\pi(1-\beta)
          A_{0}+o(1)}\int_{B_{Rr_\epsilon^{1/(1-\beta)}}}\f{e^{4\pi(1-\beta)
          u_{\epsilon,r}^2/\tau_{\epsilon,r}}}{|x|^{2\beta}}dx\\\nonumber
          &=&r^{-2(1-\beta)}e^{4\pi(1-\beta)
          A_{0}+o(1)}\int_{B_{Rr_\epsilon^{1/(1-\beta)}}}\f{e^{4\pi(1-\beta)
          u_{\epsilon,r}^2/\tau_{\epsilon,r}}-1}{|x|^{2\beta}}dx+o(1)\\\label{cc}
          &\leq&r^{-2(1-\beta)}e^{4\pi(1-\beta)
          A_{0}+o(1)}\int_{B_{r}}\f{e^{4\pi(1-\beta)
          u_{\epsilon,r}^2/\tau_{\epsilon,r}}-1}{|x|^{2\beta}}dx+o(1).
          \eea
          Combining (\ref{r2}) with (\ref{cc}), one concludes for any fixed $R>0$,
          \be\label{u-p-1}\limsup_{\epsilon\ra 0}\int_{B_{Rr_\epsilon^{1/(1-\beta)}}}\f{e^{4\pi(1-\beta-\epsilon)
          u_\epsilon^2}-1}{|x|^{2\beta}}dx\leq \f{\pi}{1-\beta}e^{1+4\pi(1-\beta)A_0}.\ee
          In view of Lemma \ref{Lemma 6}, we calculate
          \bna
          \int_{B_{Rr_\epsilon^{1/(1-\beta)}}}\f{e^{4\pi(1-\beta-\epsilon)
          u_\epsilon^2}-1}{|x|^{2\beta}}dx&=&r_\epsilon^2\int_{B_R}\f{e^{4\pi(1-\beta-\epsilon)
          u_\epsilon^2(r_\epsilon^{{1}/{(1-\beta)}}y)}}{|y|^{2\beta}}dy+o_\epsilon(1)\\
          &=&\f{\lambda_\epsilon}{c_\epsilon^2}\le(\int_{B_R}\f{e^{8\pi(1-\beta)\varphi(y)}}{|y|^{2\beta}}dy+o_\epsilon(1)\ri)+o_\epsilon(1)\\
          &=&\f{\lambda_\epsilon}{c_\epsilon^2}(1+o_R(1)+o_\epsilon(1))+o_\epsilon(1).
          \ena
          As a consequence,
          \be\label{rr-e}\lim_{R\ra\infty}\lim_{\epsilon\ra 0}\int_{B_{Rr_\epsilon^{1/(1-\beta)}}}\f{e^{4\pi(1-\beta-\epsilon)
          u_\epsilon^2}-1}{|x|^{2\beta}}dx=\lim_{\epsilon\ra 0}\f{\lambda_\epsilon}{c_\epsilon^2}.\ee
          Combining (\ref{u-p-1}), (\ref{rr-e}) and (\ref{lit}), in view of (\ref{l-t1}) and (\ref{l-t2}), we arrive at
          \be\label{upper}\sup_{u\in W^{1,2}(\mathbb{R}^2),\,\|u\|_{1,\tau}\leq 1}\int_{\mathbb{R}^2}\f{e^{4\pi(1-\beta)u^2}-1}{|x|^{2\beta}}dx=\lim_{\epsilon\ra 0}\int_{\mathbb{R}^2}\f{e^{4\pi(1-\beta-\epsilon)u_\epsilon^2}-1}{|x|^{2\beta}}dx
          \leq \f{\pi}{1-\beta}e^{1+4\pi(1-\beta)A_0}.\ee

          \subsection{Test function computation}

          We now construct test functions such that (\ref{upper}) does not hold. Precisely we construct a sequence of functions
          $\phi_\epsilon\in W^{1,2}(\mathbb{R}^2)$ satisfying $\|\phi_\epsilon\|_{1,\tau}=1$ and
          \be\label{great}\int_{\mathbb{R}^2}\f{e^{4\pi(1-\beta)\phi_\epsilon^2}-1}{|x|^{2\beta}}dx>\f{\pi}{1-\beta}
          e^{1+4\pi(1-\beta)A_0}\ee
          for sufficiently small $\epsilon>0$. For this purpose we set
          $$\phi_\epsilon(x)=\le\{
     \begin{array}{llll}
     &c+\f{1}{c}\le(-\f{1}{4\pi(1-\beta)}\log(1+\f{\pi}{1-\beta}\f{|x|^{2(1-\beta)}}
     {\epsilon^{2(1-\beta)}})+b\ri),
     \quad &x\in\overline{{B}}_{R\epsilon}\\[1.2ex]
     &\f{G}{c},\quad & x\in \mathbb{R}^2\setminus{B}_{R\epsilon},
     \end{array}
     \ri.$$
     where $G$ is given as in Lemma \ref{Green}, $R=(-\log\epsilon)^{1/(1-\beta)}$,
      $b$ and $c$ are constants depending only on $\epsilon$ to be
     determined later. To ensure $\phi_\epsilon\in W^{1,2}(\mathbb{R}^2)$, we let
     $$c+\f{1}{c}\le(-\f{1}{4\pi(1-\beta)}\log\le(1+\f{\pi}{1-\beta}R^{2(1-\beta)}\ri)+b\ri)=\f{1}{c}
     \le(-\f{1}{2\pi}\log(R\epsilon)+A_0+w(R\epsilon)\ri).$$
     This leads to
     \be\label{c2}
     c^2=\f{1}{4\pi(1-\beta)}\log\f{\pi}{1-\beta}+A_0-b-\f{1}{2\pi}\log\epsilon+O(\f{1}{R^{2-2\beta}})+O(R\epsilon).
     \ee
     Now we calculate
     \bea\nonumber
     \int_{\mathbb{R}^2\setminus B_{R\epsilon}}(|\nabla \phi_\epsilon|^2+\tau \phi_\epsilon^2)dx&=&
     \f{1}{c^2}\int_{\mathbb{R}^2\setminus B_{R\epsilon}}(|\nabla G|^2+\tau G^2)dx\\\nonumber
     &=&-\f{1}{c^2}\int_{\p B_{R\epsilon}}G\f{\p G}{\p r}d\sigma\\\label{delt}
     &=&\f{1}{c^2}\le(-\f{1}{2\pi}\log(R\epsilon)+A_0+O(R\epsilon\log(R\epsilon))\ri)
     \eea
     and
     \bea\nonumber
     \int_{B_{R\epsilon}}|\nabla\phi_\epsilon|^2dx&=&\f{1}{4\pi c^2}\int_0^{R\epsilon}
     \f{2r^{3-4\beta}}{(r^{2(1-\beta)}+\f{1-\beta}{\pi}\epsilon^{2(1-\beta)})^2}dr\\\nonumber
     &=&\f{1}{4\pi (1-\beta)c^2}\le(\log(1+\f{\pi}{1-\beta}R^{2-2\beta})+\f{1}{1+\f{\pi}{1-\beta}R^{2-2\beta}}-1\ri)\\
     &=&\f{1}{4\pi (1-\beta)c^2}\le(\log\f{\pi}{1-\beta}+\log R^{2-2\beta}-1+O(\f{1}{R^{2-2\beta}})\ri).\label{delt-1}
     \eea
     Moreover, we require $b$ to be bounded with respect to $\epsilon$. It then follows that
     \be\label{delt-2}\int_{B_{R\epsilon}}\phi_\epsilon^2dx=\f{1}{c^2}\int_{B_{R\epsilon}}
     \le(c^2-\f{1}{4\pi(1-\beta)}\log(1+\f{\pi}{1-\beta}\f{|x|^{2(1-\beta)}}
     {\epsilon^{2(1-\beta)}})+b\ri)^2dx=O(R\epsilon).\ee
     Combining (\ref{delt}), (\ref{delt-1}) and (\ref{delt-2}), we obtain
     $$
     \|\phi_\epsilon\|_{1,\tau}^2=\f{1}{c^2}\le(-\f{1}{2\pi}\log\epsilon+A_0-\f{1}{4\pi(1-\beta)}
     +\f{1}{4\pi(1-\beta)}\log\f{\pi}{1-\beta}+O(\f{1}{R^{2-2\beta}})\ri).
     $$
     Setting $\|\phi_\epsilon\|_{1,\tau}=1$, we have
     \be\label{c-squ}c^2=-\f{1}{2\pi}\log\epsilon+A_0-\f{1}{4\pi(1-\beta)}
     +\f{1}{4\pi(1-\beta)}\log\f{\pi}{1-\beta}+O(\f{1}{R^{2-2\beta}}),\ee
     which together with (\ref{c2}) leads to
     \be\label{b-2}b=\f{1}{4\pi(1-\beta)}+O(\f{1}{R^{2-2\beta}}).\ee
     For all $x\in B_{R\epsilon}$, it follows from (\ref{c-squ}) and (\ref{b-2}) that
     \bna
     4\pi(1-\beta)\phi_\epsilon^2(x)&\geq&4\pi(1-\beta)c^2+8\pi(1-\beta)b-2
     \log\le(1+\f{\pi}{1-\beta}\f{|x|^{2-2\beta}}{\epsilon^{2-2\beta}}\ri)\\
     &=&-2
     \log\le(1+\f{\pi}{1-\beta}\f{|x|^{2-2\beta}}{\epsilon^{2-2\beta}}\ri)
     -2(1-\beta)\log\epsilon\\&&
     +4\pi(1-\beta) A_0+\log\f{\pi}{1-\beta}+1+O(\f{1}{R^{2-2\beta}}).
     \ena
     Hence
     \bea\nonumber
     \int_{{B}_{R\epsilon}}\f{e^{4\pi(1-\beta)\phi_\epsilon^2}-1}{|x|^{2\beta}}dx&\geq&
     \f{\pi}{1-\beta}\epsilon^{-2(1-\beta)}e^{1+4\pi(1-\beta)A_0+O(\f{1}{R^{2-2\beta}})}\\\nonumber
     &&\times\int_{{B}_{R\epsilon}}\f{1}{(1+\f{\pi}{1-\beta}\f{|x|^{2(1-\beta)}}{\epsilon^{2(1-\beta)}})^2|x|^{2\beta}}dx
     +O((R\epsilon)^{2-2\beta})\\\nonumber
     &=&\f{\pi}{1-\beta}e^{1+4\pi(1-\beta)A_0+O(\f{1}{R^{2-2\beta}})}\\\nonumber
     &&\times\int_{{B}_{R}}\f{1}{(1+\f{\pi}{1-\beta}|y|^{2(1-\beta)})^2|y|^{2\beta}}dy+O((R\epsilon)^{2-2\beta})\\\label{inner}
     &=&\f{\pi}{1-\beta}e^{1+4\pi(1-\beta)A_0}+O(\f{1}{R^{2-2\beta}}).
     \eea
     Also we calculate
     \be\label{outer}\int_{\mathbb{R}^2\setminus B_{R\epsilon}}\f{e^{4\pi(1-\beta)\phi_\epsilon^2}-1}{|x|^{2\beta}}dx
     \geq \f{4\pi(1-\beta)}{c^2}\int_{\mathbb{R}^2\setminus B_{R\epsilon}}\f{G^2}{|x|^{2\beta}}dx=\f{4\pi(1-\beta)}{c^2}
     \le(\int_{\mathbb{R}^2}\f{G^2}{|x|^{2\beta}}dx+o_\epsilon(1)\ri).\ee
     Combining (\ref{inner}) and (\ref{outer}) and noting that $c^2/{R^{2-2\beta}}=o_\epsilon(1)$,
      we have
     $$\int_{\mathbb{R}^2}\f{e^{4\pi(1-\beta)\phi_\epsilon^2}-1}{|x|^{2\beta}}dx\geq
     \f{\pi}{1-\beta}e^{1+4\pi(1-\beta)A_0}+\f{4\pi(1-\beta)}{c^2}
     \le(\int_{\mathbb{R}^2}\f{G^2}{|x|^{2\beta}}dx+o_\epsilon(1)\ri).
     $$
     Therefore we conclude (\ref{great}) for sufficiently small $\epsilon>0$.

     \subsection{Completion of the proof of Theorem \ref{Theorem 2}}

      Comparing (\ref{great}) with (\ref{upper}), we conclude that $c_\epsilon$ must be bounded. Then applying elliptic estimates
      to (\ref{Euler-Lagrange}), we get the desired extremal function. This ends the proof of Theorem \ref{Theorem 2}. $\hfill\Box$

\section{$N$-dimensional case}\label{Sec 3}

In this section, we prove Theorems \ref{Theorem 1} and \ref{Theorem 2} in the case that $N\geq 3$. We put emphasis on
the essential difference between 2 dimensions and $N$ dimensions. In the sequel, we denote
$\Delta_Nu={\rm div}(|\nabla u|^{N-2}\nabla u)$ for any $u\in W^{1,N}(\mathbb{R}^N)$ $(N\geq 3)$.
Let $\zeta:\mathbb{N}\times\mathbb{R}\ra\mathbb{R}$ be defined as in (\ref{MF}).
Obviously, one has
\be\label{derivative}\f{d}{dt}\zeta(N,t)=\zeta(N-1,t).\ee
In view of (\cite{YangJFA}, Lemma 2.1), there holds for all $p\geq 1$ and $t\geq 0$,
\be\label{Lp}\le(\zeta(N,t)\ri)^p\leq \zeta(N,pt). \ee

\subsection{A priori estimates}

We need elliptic estimates for quasi-linear equations as below.
\begin{theorem}\label{Estimates}
Let $R>0$ be fixed. Suppose that $u\in W^{1,N}(B_R)$ is a weak solution of
$$-\Delta_N u=f\quad {\rm in}\quad B_R\subset\mathbb{R}^N.$$
Then the following a priori estimates hold:
\begin{itemize}
\item{({\bf Harnack inequality}) If $u\geq 0$ and $f\in L^p(B_R)$ for some $p>1$, then there exists some constant $C$
depending only on $N$, $R$ and $p$ such that
$\sup_{B_{R/2}}u\leq C(\inf_{B_{R/2}}u+\|f\|_{L^p(B_R)})$;}
\item{({\bf $C^\alpha$-estimate}) If $\|u\|_{L^\infty(B_R)}\leq L$ and $\|f\|_{L^p(B_R)}\leq M$ for some $p>1$, then there exists two constants $0<\alpha\leq 1$
and $C$ depending only on $N$, $R$, $p$, $L$ and $M$  such that $u\in C^\alpha(\overline{B}_{R/2})$ and
$\|u\|_{C^\alpha(\overline{B}_{R/2})}\leq C$;}
\item{({\bf $C^{1,\alpha}$-estimate}) If $\|u\|_{L^\infty(B_R)}\leq L$ and $\|f\|_{L^\infty(B_R)}\leq M$, then there exists two constants $0<\alpha\leq 1$
and $C$ depending only on $N$, $R$, $L$ and $M$  such that $u\in C^{1,\alpha}(\overline{B}_{R/2})$ and
$\|u\|_{C^{1,\alpha}(\overline{B}_{R/2})}\leq C$.}
\end{itemize}
\end{theorem}

In the above theorem, the first two estimates were obtained by J. Serrin
(\cite{Serrin}, Theorems 6 and 8), while the third estimate was proved by
Tolksdorf (\cite{Tolksdorf}, Theorem 1).

\subsection{Extremal functions for subcritical Trudinger-Moser inequalities}

In this subsection, we prove Theorem \ref{Theorem 1} in the case $N\geq 3$. The proof is based on a direct method of
variation. Throughout this section, we denote for simplicity
\be\label{beta-e}\beta_{N,\epsilon}=\alpha_N(1-\beta-\epsilon).\ee

{\it Proof of Theorem \ref{Theorem 1}.}  Let $\mathscr{S}_N$ be a subset of
$W^{1,N}(\mathbb{R}^N)$ consisting of all functions, which are nonnegative decreasing radially  symmetric almost everywhere.
By a rearrangement argument, we have
$$\Lambda_{N,\beta,\tau,\epsilon}
=\sup_{u\in \mathscr{S}_N,\,\|u\|_{1,\tau}\leq 1}\int_{\mathbb{R}^N}
\f{\zeta(N,\beta_{N,\epsilon}u^{\f{N}{N-1}})}{|x|^{N\beta}}dx,$$
where $\Lambda_{N,\beta,\tau,\epsilon}$ is defined as in (\ref{subcrit}) and $\beta_{N,\epsilon}$ is defined as in (\ref{beta-e}).
Take $u_j\in\mathscr{S}_N$ with $\|u_j\|_{1,\tau}\leq 1$ such that
$$\lim_{j\ra\infty}\int_{\mathbb{R}^N}\f{\zeta(N,\beta_{N,\epsilon}u_j^{\f{N}{N-1}})}{|x|^{N\beta}}dx
=\Lambda_{N,\beta,\tau,\epsilon}.$$
Up to a subsequence, we can find some function $u_\epsilon$ such that
$u_j$ converges to $u_\epsilon$ weakly in $W^{1,N}(\mathbb{R}^N)$, strongly in $L^q_{\rm loc}(\mathbb{R}^N)$ for any $q>0$,
and a.e. in $\mathbb{R}^N$. Obviously $u_\epsilon\in\mathscr{S}_N$. It follows from the weak convergence of $u_j$ in
$W^{1,N}(\mathbb{R}^N)$ that
\bna
\int_{\mathbb{R}^N}|\nabla u_\epsilon|^Ndx=\lim_{j\ra\infty}\int_{\mathbb{R}^N}|\nabla u_\epsilon|^{N-2}\nabla u_j\nabla u_\epsilon dx,
\ena
which together with the H\"older inequality leads to
\be\label{e-g}\int_{\mathbb{R}^N}|\nabla u_\epsilon|^Ndx\leq \limsup_{j\ra\infty}\int_{\mathbb{R}^N}|\nabla u_j|^Ndx.\ee
While it follows from $u_j\ra u_\epsilon$ in $L^q_{\rm loc}(\mathbb{R}^N)$ for any $q>0$ that for any fixed $R>0$,
\be\label{e-g1}\int_{B_R}u_\epsilon^Ndx=\lim_{j\ra\infty}\int_{B_R}u_j^Ndx.\ee
Combining (\ref{e-g}) and (\ref{e-g1}), one can easily see that
$\|u_\epsilon\|_{1,\tau}\leq \limsup_{j\ra\infty}\|u_j\|_{1,\tau}\leq 1$.
Given any $\nu>0$, there hold
\bea\label{ge-1}
\int_{|x|>\nu^{-\f{1}{N\beta}}}\f{\zeta(N,\beta_{N,\epsilon}u_j^{{N}/{(N-1)}})}{|x|^{N\beta}}dx\leq
\nu \int_{\mathbb{R}^N}\zeta(N,\beta_{N,\epsilon}u_j^{{N}/{(N-1)}})dx\leq \nu \Lambda_{N,0,\tau},\\[1.2ex]\label{ge-2}
\int_{|x|>\nu^{-\f{1}{N\beta}}}\f{\zeta(N,\beta_{N,\epsilon}u_\epsilon^{{N}/{(N-1)}})}{|x|^{N\beta}}dx\leq
\nu \int_{\mathbb{R}^N}\zeta(N,\beta_{N,\epsilon}u_\epsilon^{{N}/{(N-1)}})dx\leq  \nu\Lambda_{N,0,\tau},
\eea
where $\Lambda_{N,0,\tau}$ is defined as in (\ref{crit}).
 In view of (\ref{derivative}), we have by the mean value theorem,
\bea\nonumber\zeta(N,\beta_{N,\epsilon} u_j^{\f{N}{N-1}})-\zeta(N,\beta_{N,\epsilon} u_\epsilon^{\f{N}{N-1}})
&=&\zeta(N-1,\vartheta)\beta_{N,\epsilon}(u_j^{\f{N}{N-1}}-u_\epsilon^{\f{N}{N-1}})\\
&\leq&\max\{\zeta(N-1,\beta_{N,\epsilon} u_j^{\f{N}{N-1}}), \zeta(N-1,\beta_{N,\epsilon} u_\epsilon^{\f{N}{N-1}})\}\nonumber\\
&&\qquad\times\,
\beta_{N,\epsilon}(u_j^{\f{N}{N-1}}-u_\epsilon^{\f{N}{N-1}}),\label{mean}\eea
where $\vartheta$ lies between $\beta_{N,\epsilon} u_j^{{N}/{(N-1)}}$ and $\beta_{N,\epsilon} u_\epsilon
^{{N}/{(N-1)}}$. Employing (\ref{Lp}) and (\ref{Singular-RN}), we calculate for some $p$, $1<p<\min\{\f{1}{1-\epsilon},\f{1}{\beta}\}$,
\bna
\int_{|x|\leq \nu^{-\f{1}{N\beta}}}\f{\le(\zeta(N-1,\vartheta)\ri)^p}{|x|^{N\beta p}}dx&\leq&
\int_{|x|\leq \nu^{-\f{1}{N\beta}}}\f{\zeta(N-1,p\vartheta)}{|x|^{N\beta p}}dx\\
&=&\int_{|x|\leq \nu^{-\f{1}{N\beta}}}\f{\zeta(N,p\vartheta)}{|x|^{N\beta p}}dx+
\int_{|x|\leq \nu^{-\f{1}{N\beta}}}\f{1}{(N-2)!}\f{(p\vartheta)^{N-2}}{|x|^{N\beta p}}dx\\
&\leq&\int_{\mathbb{R}^N}\f{\zeta(N,\beta_{N,\epsilon} pu_j^{\f{N}{N-1}})}{|x|^{N\beta p}}dx+
\int_{\mathbb{R}^N}\f{\zeta(N,
\beta_{N,\epsilon} pu_\epsilon^{\f{N}{N-1}})}{|x|^{N\beta p}}dx +C_1\leq C,\ena
where $C_1$ is a constant depending only on $N,\beta, p$, while $C$ is a constant depending on $N,\beta,\epsilon$ and $p$.
This together with (\ref{ge-1})-(\ref{mean}), the H\"older inequality and
the fact that $u_j\ra u_\epsilon$ in $L^q_{\rm loc}(\mathbb{R}^N)$ for any $q>0$ implies that
$$\Lambda_{N,\beta,\tau,\epsilon}=\lim_{j\ra\infty}\int_{\mathbb{R}^N}\f{\zeta(N,\beta_{N,\epsilon} u_j^{\f{N}{N-1}})}{|x|^{N\beta}}dx=
\int_{\mathbb{R}^N}\f{\zeta(N,\beta_{N,\epsilon} u_\epsilon^{\f{N}{N-1}})}{|x|^{N\beta}}dx.$$
 Clearly we must have $\|u_\epsilon\|_{1,\tau}=1$. Moreover,
by a straightforward calculation, we derive the Euler-Lagrange equation of $u_\epsilon$ as follows:
\be\label{EL-N}\le\{\begin{array}{lll}-\Delta_N u_\epsilon+\tau u_\epsilon^{N-1}=\f{1}{\lambda_\epsilon}
\f{u_\epsilon^{{1}/{(N-1)}}}{|x|^{N\beta}}\zeta(N-1,\beta_{N,\epsilon}u_\epsilon^{{N}/{(N-1)}})\quad{\rm in}\quad \mathbb{R}^N,
\\[1.2ex]
\lambda_\epsilon=\int_{\mathbb{R}^N} |x|^{-N\beta}u_\epsilon^{{N}/{(N-1)}}\zeta(N-1,\beta_{N,\epsilon}u_\epsilon^{{N}/{(N-1)}})dx.
\end{array}\ri.\ee
Applying Theorem \ref{Estimates} to (\ref{EL-N}), we have  $u_\epsilon\in C^1(\mathbb{R}^N\setminus\{0\})\cap
C^0(\mathbb{R}^N)$. This completes the proof of the theorem. $\hfill\Box$\\

The remaining part of this section is devoted to the proof of Theorem \ref{Theorem 2}.

\subsection{Elementary properties of $u_\epsilon$}

Similar to Lemma \ref{Lemma 4}, we have the following:

\begin{lemma}\label{Lemma-N1} Let $\lambda_\epsilon$ be defined as in (\ref{EL-N}). Then there holds
$\liminf_{\epsilon\ra 0}\lambda_\epsilon>0$.
\end{lemma}

{\it Proof.} Employing the Lebesgue dominated convergence theorem and noting that $u_\epsilon$
is a maximizer for subcritical Trudinger-Moser inequalities, we have for all $u\in W^{1,N}(\mathbb{R}^N)$ with
$\|u\|_{1,\tau}\leq 1$,
\bna
\int_{\mathbb{R}^N}\f{\zeta(N,\alpha_N(1-\beta)|u|^{\f{N}{N-1}})}{|x|^{N\beta}}dx&=&
\lim_{\epsilon\ra 0}\int_{\mathbb{R}^N}\f{\zeta(N,\beta_{N,\epsilon}|u|^{\f{N}{N-1}})}{|x|^{N\beta}}dx\\
&\leq&\limsup_{\epsilon\ra0}\int_{\mathbb{R}^N}\f{\zeta(N,\beta_{N,\epsilon} u_\epsilon^{\f{N}{N-1}})}{|x|^{N\beta}}dx.
\ena
One easily concludes
\be\label{lim-N1}\Lambda_{N,\beta,\tau}=\lim_{\epsilon\ra 0}\int_{\mathbb{R}^N}
\f{\zeta(N,\beta_{N,\epsilon}u_\epsilon^{\f{N}{N-1}})}{|x|^{N\beta}}dx.\ee
Since for any $t\geq 0$,
$$t\zeta(N-1,t)=\sum_{k=N-2}^\infty\f{t^{k+1}}{k!}=\sum_{k=N-1}^\infty\f{t^k}{(k-1)!}\geq
\sum_{k=N-1}^\infty\f{t^k}{k!}=\zeta(N,t),$$
one has
\bna
\lambda_\epsilon\geq \f{1}{\beta_{N,\epsilon}}\int_{\mathbb{R}^N}\f{\zeta(N,\beta_{N,\epsilon}
u_\epsilon^{\f{N}{N-1}})}{|x|^{N\beta}}dx=
\f{1}{\alpha_N(1-\beta)}\Lambda_{N,\beta,\tau}+o_\epsilon(1).
\ena
Thus we get the desired result since $\Lambda_{N,\beta,\tau}>0$.
$\hfill\Box$\\

Since $\|u_\epsilon\|_{1,\tau}=1$, one can find some function $u_0$ such that $u_\epsilon$ converges to
$u_0$ weakly in $W^{1,N}(\mathbb{R}^N)$, strongly in $L^q_{\rm loc}(\mathbb{R}^N)$ for any $q>0$, and a.e. in
$\mathbb{R}^N$.
Denote $c_\epsilon=u_\epsilon(0)$. If $c_\epsilon$ is a bounded sequence, then applying a priori estimates in
Theorem \ref{Estimates}
to (\ref{EL-N}), we conclude that $u_\epsilon\ra u_0$ in $C^0_{\rm loc}(\mathbb{R}^N)\cap C^1_{\rm loc}(\mathbb{R}^N\setminus\{0\})$.
It is not difficult to see that
$$\int_{\mathbb{R}^N}\f{\zeta(N,\alpha_N(1-\beta)u_0^{\f{N}{N-1}})}{|x|^{N\beta}}dx=\lim_{\epsilon\ra 0}
\int_{\mathbb{R}^N}\f{\zeta(N,\alpha_N(1-\beta)u_\epsilon^{\f{N}{N-1}})}{|x|^{N\beta}}dx=\Lambda_{N,\beta,\tau}.$$
This also implies that $\|u_0\|_{1,\tau}=1$ and thus $u_0$ is the desired maximizer for the critical Trudinger-Moser functional.
In the following, without loss of generality, we assume $c_\epsilon\ra+\infty$ as $\epsilon\ra 0$.

\begin{lemma}\label{Lemma-N2}
$u_0\equiv 0$ and up to a subsequence, $|\nabla u_\epsilon|^Ndx\rightharpoonup \delta_0$ weakly in the sense of measure.
\end{lemma}

{\it Proof.} We first prove that $|\nabla u_\epsilon|^Ndx\rightharpoonup \delta_0$. Suppose not. There exists $r_0>0$ such that
$$\limsup_{\epsilon\ra 0}\int_{B_{r_0}}|\nabla u_\epsilon|^Ndx\leq \eta<1.$$
Note that $u_\epsilon$ is decreasing radially symmetric. Let $\widetilde{u}_\epsilon(x)=u_\epsilon(x)-u_\epsilon(r_0)$ for $x\in B_{r_0}$.
Then $\widetilde{u}_\epsilon\in W_0^{1,N}(B_{r_0})$ satisfies $\|\nabla \widetilde{u}_\epsilon\|_{L^N(B_{r_0})}\leq \eta<1$.
Denote
$$f_\epsilon(x)=\f{1}{\lambda_\epsilon}
\f{u_\epsilon^{{1}/{(N-1)}}(x)}{|x|^{N\beta}}\zeta(N-1,\beta_{N,\epsilon}u_\epsilon^{\f{N}{N-1}}(x)).$$
There holds for any $p>1$, $p_1>1$ and $1/p_1+1/p_2=1$,
\bea\nonumber
\int_{B_{r_0}}f_\epsilon^p(x)dx&\leq&\int_{B_{r_0}}\f{1}{\lambda_\epsilon^p}\f{u_\epsilon^{p/(N-1)}}{|x|^{N\beta p}}
\zeta(N-1,\beta_{N,\epsilon} pu_\epsilon^{\f{N}{N-1}}(x))dx\\\nonumber
&\leq&\f{1}{\lambda_\epsilon^p}\le(\int_{B_{r_0}}\f{u_\epsilon^{pp_1/(N-1)}}{|x|^{N\beta p}}dx\ri)^{1/p_1}
\le(\int_{B_{r_0}}\f{\zeta(N-1,\beta_{N,\epsilon} pp_2u_\epsilon^{\f{N}{N-1}})}{|x|^{N\beta p}}dx\ri)^{1/p_2}\\\label{r0}
&\leq&\f{1}{\lambda_\epsilon^p}\le(\int_{B_{r_0}}\f{u_\epsilon^{pp_1/(N-1)}}{|x|^{N\beta p}}dx\ri)^{1/p_1}
\le(\int_{B_{r_0}}\f{e^{\alpha_N(1-\beta) pp_2u_\epsilon^{\f{N}{N-1}}}}{|x|^{N\beta p}}dx\ri)^{1/p_2}.
\eea
Since $u_\epsilon$ is nonnegative decreasing radially symmetric, one has
$\int_{B_{r_0}}u_\epsilon^Ndx\geq u_\epsilon^N(r_0)\f{\omega_{N-1}}{N}r_0^N$. It follows that
\be\label{r0-norm}u_\epsilon(r_0)\leq \le(\f{N}{\omega_{N-1}}\ri)^{1/N}\f{\|u_\epsilon\|_{L^N(B_{r_0})}}{r_0}\leq
\le(\f{N}{\omega_{N-1}\tau}\ri)^{1/N}\f{1}{r_0}.\ee
Here we have used $\|u_\epsilon\|_{1,\tau}=1$. For any $\nu>0$, there exists some constant $C_0$ depending only on
$N$ and $\nu$ such that for all $x\in B_{r_0}$,
\be\label{ue-r0}
u_\epsilon^{\f{N}{N-1}}(x)\leq (1+\nu)\widetilde{u}_\epsilon^{\f{N}{N-1}}(x)+C_0u_\epsilon^{\f{N}{N-1}}(r_0).
\ee
Choosing $p>1$, $p_2>1$ sufficiently close to $1$ and $\nu>0$ sufficiently small such that $(1-\beta)pp_2(1+\nu)
+\beta p<1$, inserting (\ref{r0-norm}) and (\ref{ue-r0}) into (\ref{r0}),  and noting
that $u_\epsilon$ is bounded in $L^q(B_{r_0})$ for any fixed $q>0$, one can see from (\ref{Singular}) and
Lemma \ref{Lemma-N1} that $f_\epsilon$ is bounded in $L^p(B_{r_0})$. By the elliptic estimate
(Theorem \ref{Estimates}), $u_\epsilon$ is uniformly bounded in $B_{r_0/2}$ contradicting $c_\epsilon\ra+\infty$. This confirms that
$|\nabla u_\epsilon|^Ndx\rightharpoonup\delta_0$ in the sense of measure.

Next we prove $u_0\equiv 0$. It follows from $\|u_\epsilon\|_{1,\tau}=1$ and $|\nabla u_\epsilon|^Ndx\rightharpoonup\delta_0$
that $\|u_\epsilon\|_{L^N(\mathbb{R}^N)}=o_\epsilon(1)$,
which leads to
$$\int_{\mathbb{R}^N}u_0^Ndx\leq\limsup_{\epsilon\ra 0}\int_{\mathbb{R}^N}u_\epsilon^Ndx=0.$$
Therefore $u_0\equiv 0$ and the proof of the lemma is completed. $\hfill\Box$

\subsection{Blow-up analysis}

Let
\be\label{scale}r_\epsilon={\lambda_\epsilon^{1/N}}c_\epsilon^{-1/(N-1)}e^{-\beta_{N,\epsilon}c_\epsilon^{N/(N-1)}/N}.\ee
Define \be\label{p-sN}\psi_{N,\epsilon}(x)=c_\epsilon^{-1}
u_\epsilon(r_\epsilon^{{1}/{(1-\beta)}}x)\ee and
\be\label{v-p-eN}\varphi_{N,\epsilon}(x)=c_\epsilon^{{1}/{(N-1)}}
(u_\epsilon(r_\epsilon^{{1}/{(1-\beta)}}x)-c_\epsilon).\ee Analogous to Lemma \ref{Lemma 6}, we have
the following:

\begin{lemma}\label{Lemma 6N} Let $r_\epsilon$, $\psi_{N,\epsilon}$ and $\varphi_{N,\epsilon}$ be defined as in
 (\ref{scale})-(\ref{v-p-eN}). Then $(i)$ for any $\gamma<\alpha_N(1-\beta)/N$, there holds
$r_\epsilon e^{\gamma c_\epsilon^{N/(N-1)}}\ra 0$ as $\epsilon\ra 0$; $(ii)$ $\psi_{N,\epsilon}\ra 1$ in $C^1_{\rm loc}(\mathbb{R}^N
\setminus\{0\})\cap C^0_{\rm loc}(\mathbb{R}^N)$; $(iii)$
$\varphi_{N,\epsilon}\ra \varphi_N$ in
$C^1_{\rm loc}(\mathbb{R}^N\setminus\{0\})\cap C^0_{\rm loc}(\mathbb{R}^N)$, where
$$\varphi_N(x)=-\f{N-1}{\alpha_N(1-\beta)}\log\le(1+\f{\alpha_N}
{N^{N/(N-1)}(1-\beta)^{1/(N-1)}}|x|^{\f{N}{N-1}(1-\beta)}\ri).$$ Moreover
\be\label{enN1}\int_{\mathbb{R}^N}\f{e^{\f{N}{N-1}\alpha_N(1-\beta)\varphi_N}}{|x|^{N\beta}}dx=1.\ee
\end{lemma}

{\it Proof.} $(i)$
 In view of (\ref{scale}), one has
\bea\nonumber
r_\epsilon^Ne^{N\gamma c_\epsilon^{\f{N}{N-1}}}&=&c_\epsilon^{-\f{N}{N-1}}e^{-\alpha_N(1-\beta-\epsilon-\f{N\gamma}{\alpha_N})c_\epsilon^
{\f{N}{N-1}}}
\int_{\mathbb{R}^N}\f{u_\epsilon^{\f{N}{N-1}}\zeta(N-1,\beta_{N,\epsilon} u_\epsilon^{\f{N}{N-1}})}{|x|^{N\beta}}dx\\\nonumber
&\leq&
c_\epsilon^{-\f{N}{N-1}}e^{-\alpha_N(1-\beta-\epsilon-\f{N\gamma}{\alpha_N})c_\epsilon^
{\f{N}{N-1}}}
\int_{\mathbb{R}^N}\f{u_\epsilon^{\f{N}{N-1}}e^{\beta_{N,\epsilon} u_\epsilon^{\f{N}{N-1}}}}{|x|^{N\beta}}dx
\\\label{tds-N}
&\leq& c_\epsilon^{-\f{N}{N-1}}\int_{\mathbb{R}^N}\f{u_\epsilon^{\f{N}{N-1}}e^{N\gamma u_\epsilon^{\f{N}{N-1}}}}{|x|^{N\beta}}dx
\eea
Since $N\gamma<\alpha_N(1-\beta)$, one can see from (\ref{Singular}) that
$$\int_{\mathbb{R}^N}\f{u_\epsilon^{\f{N}{N-1}}e^{N\gamma u_\epsilon^{\f{N}{N-1}}}}{|x|^{N\beta}}dx\leq C$$
for some constant $C$, which together with (\ref{tds-N}) implies that
$r_\epsilon e^{\gamma c_\epsilon^{{N}/{(N-1)}}}=o_\epsilon(1)$.

$(ii)$ Clearly $\psi_{N,\epsilon}$ is a distributional solution to
\be\label{dN-psi}-\Delta_N\psi_{N,\epsilon}(x)= -\tau r_\epsilon^{\f{N}{1-\beta}}\psi_{N,\epsilon}^{N-1}(x)+c_\epsilon^{-N}|x|^{-N\beta}\psi_{N,\epsilon}^{\f{1}{N-1}}(x)
e^{-\beta_{N,\epsilon} c_\epsilon^{{N}/{(N-1)}}}\zeta(N-1,\beta_{N,\epsilon} u_\epsilon^{\f{N}{N-1}}(r_\epsilon^{\f{1}{1-\beta}}x)).\ee
 Applying Theorem \ref{Estimates}
to (\ref{dN-psi}), we have $\psi_{N,\epsilon}\ra \psi_N$ in $C^1_{\rm loc}(\mathbb{R}^N\setminus\{0\})\cap C^0_{\rm loc}(\mathbb{R}^N)$,
where $\psi_N$ is a distributional solution to $\Delta_N\psi_N=0$ in
$\mathbb{R}^N$. Clearly $\psi_N\equiv 1$ on $\mathbb{R}^N$.

$(iii)$ In view of (\ref{EL-N}), we derive the equation of $\varphi_{N,\epsilon}$ as follows.
\bea-\Delta_N\varphi_{N,\epsilon}(x)&=&g_{N,\epsilon}(x)\nonumber\\&=&
 |x|^{-N\beta}\psi_{N,\epsilon}^{\f{1}{N-1}}(x)
e^{-\beta_{N,\epsilon} c_\epsilon^{{N}/{(N-1)}}}\zeta(N-1,\beta_{N,\epsilon} u_\epsilon^{\f{N}{N-1}}(r_\epsilon^{\f{1}{1-\beta}}x))
 -\tau r_\epsilon^{\f{N}{1-\beta}}c_\epsilon^N\psi_{N,\epsilon}^{N-1}(x).\label{dN-varsi}\eea
Let $R, r$ be any two positive numbers such that $R>4r$. Clearly $g_{N,\epsilon}$ is bounded in $L^p(B_R)$ for some $p>1$.
Moreover, $-\varphi_{N,\epsilon}\geq 0$.
Theorem \ref{Estimates} implies that $\varphi_{N,\epsilon}$ is uniformly bounded in $B_{R/2}$. While
$g_{N,\epsilon}$ is bounded in $L^\infty(B_R\setminus B_r)$. Hence we have by applying Theorem \ref{Estimates} to (\ref{dN-varsi}),
$\varphi_{N,\epsilon}$ is bounded in $C^{1,\alpha}(B_{R/2}\setminus B_{2r})$ for some $0<\alpha<1$. Therefore
up to a subsequence, there exists some function $\varphi_N$ such that
$\varphi_{N,\epsilon}\ra \varphi_N$ in $C^1_{\rm loc}(\mathbb{R}^N\setminus\{0\})
\cap C^0_{\rm loc}(\mathbb{R}^N)$. To derive the equation of $\varphi_N$, we estimate
$$0\leq e^{-\beta_{N,\epsilon} c_\epsilon^{{N}/{(N-1)}}}\sum_{k=0}^{N-3}\f{\beta_{N,\epsilon}^k u_\epsilon^{\f{Nk}
{N-1}}(r_\epsilon^{\f{1}{1-\beta}}x)}{k!}\leq e^{-\beta_{N,\epsilon} c_\epsilon^{{N}/{(N-1)}}}\sum_{k=0}^{N-3}
\f{\beta_{N,\epsilon}^k c_\epsilon^{\f{Nk}
{N-1}}}{k!}=o_\epsilon(1)$$
uniformly on $B_R$ for any $R>0$. Moreover, by the mean value theorem, we have
\bea\nonumber
u_\epsilon^{\f{N}{N-1}}(r_\epsilon^{\f{1}{1-\beta}}x)- c_\epsilon^{\f{N}{N-1}}&=&\f{N}{N-1}\xi_\epsilon^{\f{1}{N-1}}
(u_\epsilon(r_\epsilon^{\f{1}{1-\beta}}x)- c_\epsilon)\\&=&\f{N}{N-1}(\xi_\epsilon/c_\epsilon)^{\f{1}{N-1}}\varphi_{N,\epsilon}(x)
\nonumber\\&=&\f{N}{N-1}\varphi_N(x)+o_\epsilon(1),\label{bNN}\eea
where $\xi_\epsilon$ lies between $u_\epsilon(r_\epsilon^{{1}/{(1-\beta)}}x)$ and $c_\epsilon$, and $o_\epsilon(1)\ra 0$ uniformly
on $B_R$ for any fixed $R>0$. Hence
\bna
e^{-\beta_{N,\epsilon} c_\epsilon^{{N}/{(N-1)}}}\zeta(N-1,\beta_{N,\epsilon} u_\epsilon^{\f{N}{N-1}}(r_\epsilon^{\f{1}{1-\beta}}x))
=e^{\alpha_N(1-\beta)\f{N}{N-1}\varphi_N(x)}+o_\epsilon(1).
\ena
Furthermore, we obtain the equation of $\varphi_N$
as follows:
\be\label{bubb-N}\le\{\begin{array}{lll}-\Delta_N\varphi_N(x)=\f{e^{\alpha_N(1-\beta)\f{N}{N-1}\varphi_N(x)}}
{|x|^{N\beta}}\quad{\rm in}\quad \mathbb{R}^N,\\[1.5ex]
\varphi_N(0)=\max_{\mathbb{R}^N}\varphi_N=0.\end{array}
\ri.\ee
Since $\varphi_{N,\epsilon}$ is decreasingly symmetric on $\mathbb{R}^N$, $\varphi_N$ is also decreasingly symmetric.
Denote $\varphi_N(r)=\varphi_N(x)$, where $r=|x|$ and $x\in \mathbb{R}^N$.
Then (\ref{bubb-N}) can be  reduced to an ordinary differential equation, namely
\be\label{bubb-N-r}\le\{\begin{array}{lll}\le((-r\varphi_N^\prime(r))^{N-1}\ri)^\prime=r^{N-1-N\beta}
e^{\alpha_N(1-\beta)\f{N}{N-1}\varphi_N(r)}\\[1.5ex]
\varphi_N(0)=0.\end{array}
\ri.\ee
By a standard uniqueness result of ordinary differential equations (see for example \cite{Li-Ruf}),
we can solve (\ref{bubb-N-r}) as
$$\varphi_N(r)=-\f{N-1}{\alpha_N(1-\beta)}\log\le(1+c_Nr^{\f{N}{N-1}(1-\beta)}\ri),$$
where $c_N=\alpha_NN^{-{N}/{(N-1)}}(1-\beta)^{-{1}/{(N-1)}}$.
It then follows that
\bea
\int_{\mathbb{R}^N}\f{e^{\alpha_N(1-\beta)\f{N}{N-1}\varphi_N(x)}}{|x|^{N\beta}}dx&=&
\omega_{N-1}\int_0^\infty\f{r^{N-1-N\beta}}{(1+c_Nr^{N(1-\beta)/(N-1)})^N}
dr\nonumber\\&=&\omega_{N-1}\f{N-1}{N(1-\beta)}\int_0^\infty\f{t^{N-2}}{(1+c_Nt)^N}dt.\label{b-energy}
\eea
Integration by parts gives
\bna
I_N&=&(N-1)\int_0^\infty\f{t^{N-2}}{(1+c_Nt)^N}dt\\
&=&-\f{1}{c_N}\int_0^\infty t^{N-2}d(1+c_Nt)^{1-N}\\
&=&\le.-\f{1}{c_N}t^{N-2}(1+c_Nt)^{1-N}\ri|_0^\infty+\f{N-2}{c_N}\int_0^\infty t^{N-3}(1+c_Nt)^{1-N}dt\\
&=&\f{1}{c_N}I_{N-1}.
\ena
Iteration leads to
\be\label{IN}
I_N=\f{1}{c_N^{N-2}}I_2
=\f{1}{c_N^{N-2}}\int_0^\infty\f{1}{(1+c_Nt)^2}dt=\f{1}{c_N^{N-1}}.
\ee
Inserting (\ref{IN}) into (\ref{b-energy}), we obtain
$$\int_{\mathbb{R}^N}\f{e^{\f{N}{N-1}\alpha_N(1-\beta)\varphi_N(x)}}{|x|^{N\beta}}dx=
\f{\omega_{N-1}}{N(1-\beta)}\f{1}{c_N^{N-1}}=1.$$
This completes the proof of the lemma.
 $\hfill\Box$\\

For any $0<\gamma<1$, we set $u_{\epsilon,\gamma}=\min\{u_\epsilon,\gamma c_\epsilon\}$. Then we have the following:

\begin{lemma}\label{Lemma 18}
For any $0<\gamma<1$, there holds $\lim_{\epsilon\ra 0}\int_{\mathbb{R}^N}|\nabla u_{\epsilon,\gamma}|^Ndx=\gamma$.
\end{lemma}

{\it Proof.} Testing the equation (\ref{EL-N})
by $u_{\epsilon,\gamma}$, we have for any fixed $R>0$,
\bna
\int_{\mathbb{R}^N}|\nabla u_{\epsilon,\gamma}|^Ndx&=&-\tau\int_{\mathbb{R}^N}u_\epsilon^{N-1} u_{\epsilon,\gamma}dx+
\f{1}{\lambda_\epsilon}\int_{\mathbb{R}^N} u_{\epsilon,\gamma}\f{u_\epsilon^{1/(N-1)}}{|x|^{N\beta}}
\zeta(N-1,\beta_{N,\epsilon} u_\epsilon^{N/(N-1)})dx\\
&\geq&\f{1}{\lambda_\epsilon}\int_{B_{Rr_\epsilon^{1/(1-\beta)}}}\gamma c_\epsilon \f{u_\epsilon^{1/(N-1)}}{|x|^{N\beta}}
e^{\beta_{N,\epsilon} u_\epsilon^{N/(N-1)}}dx+o_\epsilon(1)\\
&=&(1+o_\epsilon(1))\gamma\int_{B_R}\f{e^{\alpha_N(1-\beta)\f{N}{N-1}\varphi_N}}{|x|^{N\beta}}dx+o_\epsilon(1).
\ena
Hence
$$\liminf_{\epsilon\ra 0}\int_{\mathbb{R}^N}|\nabla u_{\epsilon,\gamma}|^Ndx\geq \gamma\int_{B_R}
\f{e^{\alpha_N(1-\beta)\f{N}{N-1}\varphi_N}}{|x|^{N\beta}}dx.$$
In view of (\ref{enN1}), passing to the limit $R\ra+\infty$, we obtain
\be\label{geq-N}\liminf_{\epsilon\ra 0}\int_{\mathbb{R}^N}|\nabla u_{\epsilon,\gamma}|^Ndx\geq \gamma.\ee
Similarly we have
\be\label{leq-N}\liminf_{\epsilon\ra 0}\int_{\mathbb{R}^N}|\nabla (u_\epsilon-\gamma c_\epsilon)^+|^Ndx\geq 1-\gamma.\ee
Noting that $\|u_\epsilon\|_{L^N(\mathbb{R}^N)}=o_\epsilon(1)$, we have
\be\label{1-N}\int_{\mathbb{R}^N}|\nabla u_{\epsilon,\gamma}|^Ndx+
\int_{\mathbb{R}^N}|\nabla (u_\epsilon-\gamma c_\epsilon)^+|^Ndx=\int_{\mathbb{R}^N}|\nabla u_{\epsilon}|^Ndx=1+o_\epsilon(1).\ee
Combining (\ref{geq-N})-(\ref{1-N}), we conclude the lemma. $\hfill\Box$

\begin{lemma}\label{Lemma 19} We have
\be\label{limsup-N}
\lim_{\epsilon\ra 0}\int_{\mathbb{R}^N}\f{\zeta(N,\beta_{N,\epsilon} u_\epsilon^{N/(N-1)})}{|x|^{N\beta}}dx
=\lim_{\epsilon\ra 0}\f{\lambda_\epsilon}{c_\epsilon^{N/(N-1)}}.
\ee
As a consequence, for any $\theta<N/(N-1)$, there holds $\lambda_\epsilon/c_\epsilon^\theta\ra +\infty$ as
$\epsilon\ra 0$.
\end{lemma}

{\it Proof.}
Let $0<\gamma<1$ be fixed and $u_{\epsilon,\gamma}$ be defined as before.
Applying the mean value theorem to the function $\zeta(N,t)$ and recalling (\ref{derivative}), we have
$$\zeta(N,\beta_{N,\epsilon} u_{\epsilon,\gamma}^{N/(N-1)})=\zeta(N-1,\xi_\epsilon)\beta_{N,\epsilon} u_{\epsilon,\gamma}^{N/(N-1)}
\leq \zeta(N-1,\beta_{N,\epsilon} u_{\epsilon,\gamma}^{N/(N-1)})\beta_{N,\epsilon} u_{\epsilon,\gamma}^{N/(N-1)},$$
where $\xi_\epsilon$ lies between $\beta_{N,\epsilon} u_{\epsilon,\gamma}^{N/(N-1)}$ and $0$. Since $\zeta(N-1,t)=\zeta(N,t)+t^{N-2}/(N-2)!$
for all $t\geq 0$, it follows from the above inequality that
\be\label{ineq}\zeta(N,\beta_{N,\epsilon} u_{\epsilon,\gamma}^{N/(N-1)})\leq \zeta(N,\beta_{N,\epsilon} u_{\epsilon,\gamma}^{N/(N-1)})
\beta_{N,\epsilon} u_{\epsilon,\gamma}^{N/(N-1)}+\beta_{N,\epsilon}^{N-1}u_\epsilon^N/(N-2)!.\ee
It is easy to see that
\be\label{zer}\int_{\mathbb{R}^N}\f{u_\epsilon^q}{|x|^{N\beta}}dx=o_\epsilon(1),\quad\forall q\geq N.\ee
In view of Lemma \ref{Lemma 18}, one can find some $p>1$ such that
\be\label{gamma}\limsup_{\epsilon\ra 0}\int_{\mathbb{R}^N}\f{\zeta(N,p\beta_{N,\epsilon} u_{\epsilon,\gamma}^{N/(N-1)})}
{|x|^{N\beta}}dx<\infty.\ee
By the H\"older inequality and (\ref{Lp}), one has
\be\label{hold}\int_{\mathbb{R}^N}\f{\zeta(N,\beta_{N,\epsilon} u_{\epsilon,\gamma}^{\f{N}{N-1}})
\beta_{N,\epsilon} u_{\epsilon,\gamma}^{\f{N}{N-1}}}{|x|^{N\beta}}dx\leq\le(\int_{\mathbb{R}^N}\f{\zeta(N,p\beta_{N,\epsilon}
u_{\epsilon,\gamma}^{\f{N}{N-1}})}
{|x|^{N\beta}}dx\ri)^{{1}/{p}}\le(\int_{\mathbb{R}^N}\f{(\beta_{N,\epsilon} u_{\epsilon,\gamma}^{\f{N}{N-1}})^{p^\prime}}
{|x|^{N\beta}}dx\ri)^{{1}/{p^\prime}},\ee
where $1/p+1/p^\prime=1$.
Combining (\ref{ineq})-(\ref{hold}), one concludes
\be\label{ts-0}\lim_{\epsilon\ra 0}\int_{\mathbb{R}^N}\f{\zeta(N,\beta_{N,\epsilon} u_{\epsilon,\gamma}^{\f{N}{N-1}})}
{|x|^{N\beta}}dx=0.\ee

Since $u_\epsilon\ra 0$ in $L^q_{\rm loc}(\mathbb{R}^N)$ for any $q>0$, we obtain
\bea\nonumber
\int_{u_\epsilon> \gamma c_\epsilon}\f{\zeta(N,\beta_{N,\epsilon} u_\epsilon^{\f{N}{N-1}})}{|x|^{N\beta}}dx&=&
\int_{u_\epsilon> \gamma c_\epsilon}\f{{e^{\beta_{N,\epsilon} u_\epsilon^{\f{N}{N-1}}}}}{|x|^{N\beta}}dx+o_\epsilon(1)\\
\nonumber&\leq&\f{1}{\gamma^{\f{N}{N-1}}}\int_{u_\epsilon> \gamma c_\epsilon}\f{u_\epsilon^{\f{N}{N-1}}}{c_\epsilon^{\f{N}{N-1}}}
\f{{e^{\beta_{N,\epsilon} u_\epsilon^{\f{N}{N-1}}}}}{|x|^{N\beta}}dx+o_\epsilon(1)
\\\label{gam2N}
&\leq& \f{1}{\gamma^{\f{N}{N-1}}}\f{\lambda_\epsilon}{c_\epsilon^{\f{N}{N-1}}}+o_\epsilon(1).
\eea
Combining (\ref{ts-0}) and (\ref{gam2N}), we have
$$\lim_{\epsilon\ra 0}\int_{\mathbb{R}^N}\f{\zeta(N,\beta_{N,\epsilon} u_\epsilon^{\f{N}{N-1}})}{|x|^{N\beta}}dx\leq
\f{1}{\gamma^{N/(N-1)}}\liminf_{\epsilon\ra 0}\f{\lambda_\epsilon}{c_\epsilon^{N/(N-1)}}.$$
Letting $\gamma\ra 1$, we conclude
\be\label{ll-N}
\lim_{\epsilon\ra 0}\int_{\mathbb{R}^N}\f{\zeta(N,\beta_{N,\epsilon} u_\epsilon^{\f{N}{N-1}})}{|x|^{N\beta}}dx\leq
\liminf_{\epsilon\ra 0}\f{\lambda_\epsilon}{c_\epsilon^{N/(N-1)}}.
\ee
An obvious analog of (\ref{gam3}) is
\be\label{gam3-N}\limsup_{\epsilon\ra 0}\f{\lambda_\epsilon}{c_\epsilon^{N/(N-1)}}\leq
\lim_{\epsilon\ra 0}\int_{\mathbb{R}^N}\f{\zeta(N,\beta_{N,\epsilon} u_\epsilon^{\f{N}{N-1}})}{|x|^{N\beta}}dx.\ee
Combining (\ref{ll-N}) and (\ref{gam3-N}), we obtain (\ref{limsup-N}), which together with (\ref{lim-N1})
implies that ${\lambda_\epsilon}/{c_\epsilon^{N/(N-1)}}$ has a positive lower bound. Then
 for any $\theta<N/(N-1)$, there holds $$\lambda_\epsilon/c_\epsilon^{\theta}=c_\epsilon^{N/(N-1)-\theta}
 \lambda_\epsilon/c_\epsilon^{N/(N-1)}\ra\infty.$$
 This proves the second assertion of the lemma. $\hfill\Box$

\begin{lemma}\label{Lemma 20}
$c_\epsilon^{\f{1}{N-1}}u_\epsilon\ra G$ in $C^1_{\rm loc}(\mathbb{R}^N\setminus\{0\})$ and weakly in $W^{1,q}(\mathbb{R}^N)$
for any $1<q<N$, where $G$ is a distributional solution to
$$-\Delta_NG+\tau G^{N-1}=\delta_0\quad{\rm in}\quad \mathbb{R}^N.$$
Moreover, $G\in W^{1,N}(\mathbb{R}^{N}\setminus B_r)$ for any $r>0$ and $G$ takes the form
$$G(x)=-\f{N}{\alpha_N}\log |x|+A_0+w(x),$$
where $A_0$ is a constant, and $w\in C^0(\mathbb{R}^N)\cap C^1(\mathbb{R}^N\setminus\{0\})$
satisfies $w(x)=O(|x|^N\log^{N-1}|x|)$ as $|x|\ra 0$.
\end{lemma}

{\it Proof.}
Multiplying both sides of (\ref{EL-N}) by $c_\epsilon$, we have
$$
-\Delta_N(c_\epsilon^{\f{1}{N-1}}u_\epsilon)+\tau c_\epsilon u_\epsilon^{N-1}=
\f{c_\epsilon u_\epsilon^{\f{1}{N-1}}}{\lambda_\epsilon}\f{\zeta(N-1,\beta_{N,\epsilon}
u_\epsilon^{\f{N}{N-1}})}{|x|^{N\beta}}\quad{\rm in}\quad \mathbb{R}^N.
$$
Replacing Lemma \ref{Lemma 6} and Corollary \ref{Cor} with Lemma \ref{Lemma 6N} and Lemma \ref{Lemma 19}
 respectively
in the proof of Lemma \ref{Lemma 9}, we obtain for any $\phi\in C^1_0(\mathbb{R}^N)$,
$$
\lim_{\epsilon\ra 0}\int_{\mathbb{R}^N}\f{c_\epsilon u_\epsilon^{\f{1}{N-1}}}{\lambda_\epsilon}
\f{\zeta(N-1,\beta_{N,\epsilon} u_\epsilon^{\f{N}{N-1}})}{|x|^{N\beta}}\phi dx=\phi(0).
$$
Since the remaining part of the proof is completely analogous to that of
(\cite{Li-Ruf}, Proposition 3.7 and Lemma 3.8), we omit the details but refer the reader to \cite{Li-Ruf}. $\hfill\Box$\\

To estimate the supremum $\Lambda_{N,\beta,\tau}$, we need the following:

\begin{lemma}{\label{Lemma 21}}
Let $w_\epsilon\in W_0^{1,N}(B_r)$ satisfy $\int_{B_r}|\nabla w_\epsilon|^Ndx\leq 1$, $w_\epsilon\rightharpoonup 0$
weakly in $W_0^{1,N}(B_r)$, and $w_\epsilon$ is nonnegative and radially symmetric. Then
\be\label{r-N}\limsup_{\epsilon\ra 0}\int_{B_r}\f{e^{\alpha_N(1-\beta)w_\epsilon^{{N}/{(N-1)}}}-1}{|x|^{N\beta}}dx\leq
\f{1}{1-\beta}\f{\omega_{N-1}}{N}r^{N(1-\beta)}e^{\sum_{k=1}^{N-1}\f{1}{k}}.\ee
\end{lemma}

{\it Proof.} We first prove (\ref{r-N}) for $r=1$.

Denote $w_\epsilon(|x|)=w_\epsilon(x)$. Let $v_\epsilon(x)=(1-\beta)^{{(N-1)}/{N}}w_\epsilon(|x|^{{1}/{(1-\beta)}})$. Then
$$\int_{B_1}|\nabla v_\epsilon|^Ndx=\int_{B_1}|\nabla w_\epsilon|^Ndx.$$
Clearly we can assume up to a subsequence, $v_\epsilon\rightharpoonup v_0$ weakly in $W_0^{1,N}(B_1)$, $v_\epsilon\ra v_0$
strongly in $L^N(B_1)$, and $v_\epsilon\ra v_0$ a.e. in $B_1$. Also, we can assume $w_\epsilon\ra 0$ a.e. in $B_1$. Hence
we conclude $v_0=0$ a.e. in $B_1$. By a change of variable $t=s^{1/(1-\beta)}$, there holds
\bna
\int_{B_1}\f{e^{\alpha_N(1-\beta)w_\epsilon^{\f{N}{N-1}}}-1}{|x|^{N\beta}}dx&=&\int_0^1\f{e^{\alpha_N(1-\beta)w_\epsilon^{\f{N}{N-1}}(t)}-1}
{t^{N\beta}}\omega_{N-1}t^{N-1}dt\\
&=&\f{1}{1-\beta}\int_0^1 (e^{\alpha_N(1-\beta)w_\epsilon^{\f{N}{N-1}}(s^{1/(1-\beta)})}-1)\omega_{N-1}s^{N-1}ds\\
&=&\f{1}{1-\beta}\int_0^1 (e^{\alpha_Nv_\epsilon^{\f{N}{N-1}}(s)}-1)\omega_{N-1}s^{N-1}ds\\
&=&\f{1}{1-\beta}\int_{B_1}(e^{\alpha_N v_\epsilon^{\f{N}{N-1}}}-1)dx.
\ena
This together with Lemma \ref{C-C-Lemma}  implies that
\be\label{bound-N}\limsup_{\epsilon\ra 0}\int_{B_1}\f{e^{\alpha_N(1-\beta)w_\epsilon^{\f{N}{N-1}}}-1}{|x|^{N\beta}}dx
\leq \f{1}{1-\beta}\f{\omega_{N-1}}{N}e^{\sum_{k=1}^{N-1}\f{1}{k}}.\ee

We next prove (\ref{r-N}) for the case of general $r$. Set $\tilde{w}_\epsilon(x)=w_\epsilon(rx)$ for $x\in B_1$.
One can check that
$$\int_{B_1}|\nabla \tilde{w}_\epsilon|^Ndx=\int_{B_r}|\nabla w_\epsilon|^Ndx$$
and that
\bna
\int_{B_r}\f{e^{\alpha_N(1-\beta)w_\epsilon^{\f{N}{N-1}}}-1}{|x|^{N\beta}}dx=r^{N(1-\beta)}\int_{B_1}
\f{e^{\alpha_N(1-\beta)\tilde{w}_\epsilon^{\f{N}{N-1}}}-1}
{|x|^{N\beta}}dx.
\ena
This together with (\ref{bound-N}) gives the desired result. $\hfill\Box$\\

By the equation (\ref{EL-N}) and $\|u_\epsilon\|_{1,\tau}=1$, we have
\bea\nonumber
\int_{B_r}|\nabla u_\epsilon|^Ndx&=&1-\int_{\mathbb{R}^N\setminus B_r}(|\nabla u_\epsilon|^N+\tau u_\epsilon^N)dx
-\tau\int_{B_r}u_\epsilon^Ndx\\\nonumber
&=&1-\int_{\mathbb{R}^N\setminus B_r}\f{u_\epsilon^{\f{N}{N-1}}}{\lambda_\epsilon}
\f{\zeta(N-1,\beta_{N,\epsilon}u_\epsilon^{\f{N}{N-1}})}{|x|^{N\beta}}dx\\\label{shang-N}
&&
+\int_{\p B_r}u_\epsilon|\nabla u_\epsilon|^{N-2}\f{\p u_\epsilon}{\p r}d\sigma-\tau\int_{B_r}u_\epsilon^Ndx.
\eea
We estimate the right three terms on the above equation respectively.
The first term can be calculated by
\bea\nonumber
\int_{\mathbb{R}^N\setminus B_r}\f{u_\epsilon^{\f{N}{N-1}}}{\lambda_\epsilon}\f{\zeta(N-1,\beta_{N,\epsilon}u_\epsilon^{\f{N}{N-1}})}{|x|^{N\beta}}dx
&=&
\f{1}{c_\epsilon^{N/(N-1)}}\f{c_\epsilon^{N/(N-1)}}{\lambda_\epsilon}\int_{\mathbb{R}^N\setminus B_r}{u_\epsilon^{\f{N}{N-1}}}
\f{\zeta(N-1,\beta_{N,\epsilon}u_\epsilon^{\f{N}{N-1}})}{|x|^{N\beta}}dx\\
&=&\f{o_\epsilon(1)}{c_\epsilon^{N/(N-1)}}.\label{First}
\eea
A straightforward calculation on the second term reads
\bea\int_{\p B_r}u_\epsilon|\nabla u_\epsilon|^{N-2}\f{\p u_\epsilon}{\p r}d\sigma
&=&\f{1}{c_\epsilon^{N/(N-1)}}\le(\int_{\p B_r}G|\nabla G|^{N-2}\f{\p G}{\p r}d\sigma+o_\epsilon(1)\ri)\nonumber\\
&=&\f{1}{c_\epsilon^{N/(N-1)}}\le(G(r)\int_{B_r}\Delta_N Gdx+o_\epsilon(1)\ri)\nonumber\\\label{Second}
&=&\f{1}{c_\epsilon^{N/(N-1)}}\le(-G(r)+\tau G(r)\int_{B_r}G^{N-1}dx+o_\epsilon(1)\ri),\eea
since $G$ is a distributional solution of $-\Delta_NG+\tau G^{N-1}=\delta_0$. Concerning the third term, one has
\be\label{Third}\int_{B_r}u_\epsilon^Ndx=\f{1}{c_\epsilon^{{N}/{(N-1)}}}\le(\int_{B_r}G^Ndx+o_\epsilon(1)\ri).\ee
Inserting (\ref{First})-(\ref{Third}) into (\ref{shang-N}) and noting that $G(x)=-\f{N}{\alpha_N}\log|x|+A_0+w(x)$, we conclude
\be\label{r-energy-N}\int_{B_r}|\nabla u_\epsilon|^Ndx=1-\f{1}{c_\epsilon^{N/(N-1)}}\le(\f{N}{\alpha_N}\log \f{1}{r}+A_0+o_\epsilon(1)+o_r(1)\ri).\ee

Define $u_{\epsilon,r}=(u_\epsilon-u_\epsilon(r))^+$, the positive part of $u_\epsilon-u_\epsilon(r)$. Obviously
$u_{\epsilon,r}\in W_0^{1,N}(B_r)$. It follows from Lemma \ref{Lemma 21} that
\be\label{llqN}
\limsup_{\epsilon\ra 0}\int_{B_r}\f{e^{\alpha_N(1-\beta)u_{\epsilon,r}^{{N}/{(N-1)}}/\tau_{\epsilon,r}}-1}{|x|^{N\beta}}dx
\leq \f{1}{1-\beta}\f{\omega_{N-1}}{N}r^{N(1-\beta)}e^{\sum_{k=1}^{N-1}\f{1}{k}},
\ee
where $\tau_{\epsilon,r}=\|\nabla u_\epsilon\|_{L^N(B_r)}^{N/(N-1)}$. One can see from Lemma \ref{Lemma 6N} that
$u_\epsilon=c_\epsilon+o_\epsilon(1)$ on $B_{Rr_\epsilon^{1/(1-\beta)}}$. This together with Lemma \ref{Lemma 20} and (\ref{r-energy-N}) leads to
that on $B_{Rr_\epsilon^{1/(1-\beta)}}\subset B_r$,
\bna
\beta_{N,\epsilon} u_\epsilon^{\f{N}{N-1}}&\leq&\alpha_N(1-\beta)(u_{\epsilon,r}+u_\epsilon(r))^{\f{N}{N-1}}\\
&=&\alpha_N(1-\beta)u_{\epsilon,r}^{\f{N}{N-1}}+\f{N}{N-1}\alpha_N(1-\beta)u_{\epsilon,r}^{\f{1}{N-1}}u_\epsilon(r)+o_\epsilon(1)\\
&=&\alpha_N(1-\beta)u_{\epsilon,r}^{\f{N}{N-1}}+\f{N}{N-1}\alpha_N(1-\beta)G(r)+o_\epsilon(1)\\
&=&\alpha_N(1-\beta)u_{\epsilon,r}^{\f{N}{N-1}}+\f{N}{N-1}\alpha_N(1-\beta)\le(\f{N}{\alpha_N}\log\f{1}{r}+A_0\ri)
+o_r(1)+o_\epsilon(1)\\
&=&\alpha_N(1-\beta)u_{\epsilon,r}^{\f{N}{N-1}}/\tau_{\epsilon,r}+N(1-\beta)\log\f{1}{r}+\alpha_N(1-\beta)A_0+o_r(1)+o_\epsilon(1).
\ena
This together with (\ref{llqN}) leads to
\bea\nonumber
\int_{B_{Rr_\epsilon^{1/(1-\beta)}}}\f{e^{\beta_{N,\epsilon} u_\epsilon^{\f{N}{N-1}}}-1}{|x|^{N\beta}}dx&\leq&
r^{-N(1-\beta)}e^{\alpha_N(1-\beta)A_0+o(1)}\int_{B_{Rr_\epsilon^{1/(1-\beta)}}}
\f{e^{\alpha_N(1-\beta) u_{\epsilon,r}^{\f{N}{N-1}}/\tau_{\epsilon,r}}}{|x|^{N\beta}}dx\\\nonumber
&=&r^{-N(1-\beta)}e^{\alpha_N(1-\beta)A_0+o(1)}\int_{B_{Rr_\epsilon^{1/(1-\beta)}}}
\f{e^{\alpha_N(1-\beta) u_{\epsilon,r}^{\f{N}{N-1}}/\tau_{\epsilon,r}}-1}{|x|^{N\beta}}dx+o(1)\\\label{105}
&\leq&\f{1}{1-\beta}\f{\omega_{N-1}}{N}e^{\sum_{k=1}^{N-1}\f{1}{k}+\alpha_N(1-\beta)A_0}+o(1).
\eea

In view of (\ref{bNN}), we obtain
\bna
\int_{B_{Rr_\epsilon^{1/(1-\beta)}}}\f{\zeta(N,\beta_{N,\epsilon} u_\epsilon^{\f{N}{N-1}})}{|x|^{N\beta}}dx
&=&r_\epsilon^N\int_{B_R}\f{e^{\beta_{N,\epsilon} u_\epsilon^{\f{N}{N-1}}(r_\epsilon^{\f{1}{1-\beta}}y)}}{|y|^{N\beta}}dy+o_\epsilon(1)\\
&=&\f{\lambda_\epsilon}{c_\epsilon^{N/(N-1)}}\le(\int_{B_R}\f{e^{\alpha_N(1-\beta)\f{N}{N-1}\varphi_N(y)}}{|y|^{N\beta}}dy
+o_\epsilon(1)\ri)+o_\epsilon(1)\\
&=&\f{\lambda_\epsilon}{c_\epsilon^{N/(N-1)}}(1+o_R(1)+o_\epsilon(1))+o_\epsilon(1).
\ena
Therefore
\be\label{limRe}\lim_{R\ra\infty}\lim_{\epsilon\ra 0}\int_{B_{Rr_\epsilon^{1/(1-\beta)}}}
\f{\zeta(N,\beta_{N,\epsilon} u_\epsilon^{\f{N}{N-1}})}{|x|^{N\beta}}dx=\lim_{\epsilon\ra 0}\f{\lambda_\epsilon}{c_\epsilon^{N/(N-1)}}.\ee
Combining (\ref{105}), (\ref{limRe}) and (\ref{limsup-N}), we conclude
\be\label{u-b-N}\Lambda_{N,\beta,\tau}=\lim_{\epsilon\ra 0}\int_{\mathbb{R}^N}\f{\zeta(N,\beta_{N,\epsilon} u_\epsilon^{\f{N}{N-1}})}{|x|^{N\beta}}dx\leq
\f{1}{1-\beta}\f{\omega_{N-1}}{N}e^{\sum_{k=1}^{N-1}\f{1}{k}+\alpha_N(1-\beta)A_0}.\ee

\subsection{Test function computation}

We now construct test functions such that (\ref{u-b-N}) does not hold. Precisely we construct a sequence of functions
          $\phi_\epsilon\in W^{1,N}(\mathbb{R}^N)$ satisfying $\|\phi_\epsilon\|_{1,\tau}=1$ and
          \be\label{great-N}\int_{\mathbb{R}^N}\f{\zeta(N,\alpha_N(1-\beta) \phi_\epsilon^{\f{N}{N-1}})}{|x|^{N\beta}}dx>
\f{1}{1-\beta}\f{\omega_{N-1}}{N}e^{\sum_{k=1}^{N-1}\f{1}{k}+\alpha_N(1-\beta)A_0}\ee
          for sufficiently small $\epsilon>0$. For this purpose we set
          $$\phi_\epsilon(x)=\le\{
     \begin{array}{llll}
     &c+\f{1}{c^{1/(N-1)}}\le(-\f{N-1}{\alpha_N(1-\beta)}\log(1+c_{N}({|x|}/{\epsilon})^{{\f{N}{N-1}(1-\beta)}})
     +b\ri),
     \quad &x\in\overline{{B}}_{R\epsilon}\\[1.2ex]
     &\f{G}{c^{1/(N-1)}},\quad & x\in \mathbb{R}^N\setminus{B}_{R\epsilon},
     \end{array}
     \ri.$$
     where $c_{N}=\alpha_N/(N^{N/(N-1)}(1-\beta)^{1/(N-1)})$, $G$ is given as in Lemma \ref{Lemma 20}, $R=(-\log\epsilon)^{1/(1-\beta)}$,
      $b$ and $c$ are constants depending only on $\epsilon$ and $\beta$ to be
     determined later. Note that $G\in W^{1,N}(\mathbb{R}^N\setminus B_r)$ for any $r>0$.
     To ensure $\phi_\epsilon\in W^{1,N}(\mathbb{R}^N)$, we let
     $$c+\f{1}{c^{1/(N-1)}}\le(-\f{N-1}{\alpha_N(1-\beta)}\log(1+c_{N}R^{{\f{N}{N-1}(1-\beta)}})
     +b\ri)=\f{G(R\epsilon)}{c^{1/(N-1)}}.$$
     By Lemma \ref{Lemma 20}, we have $G(x)=-(N/\alpha_N)\log|x|+A_0+w(x)$, where $w(x)=O(|x|^N\log^{N-1}|x|)$ as
     $|x|\ra 0$.
     Then the above equality leads to
     \be\label{cN}
     c^{\f{N}{N-1}}=\f{1}{\alpha_N(1-\beta)}\log\f{\omega_{N-1}}{N(1-\beta)}
     +A_0-b-\f{N}{\alpha_N}\log\epsilon+O(R^{-\f{N}{N-1}(1-\beta)}).
     \ee
     Now we calculate by the equation of $G$,
     \bea\nonumber
     \int_{\mathbb{R}^N\setminus B_{R\epsilon}}(|\nabla \phi_\epsilon|^N+\tau \phi_\epsilon^N)dx&=&
     \f{1}{c^{\f{N}{N-1}}}\int_{\mathbb{R}^N\setminus B_{R\epsilon}}(|\nabla G|^N+\tau G^N)dx\\\nonumber
     &=&-\f{1}{c^{\f{N}{N-1}}}\int_{\p B_{R\epsilon}}G|\nabla G|^{N-2}\f{\p G}{\p\nu}d\sigma\\\nonumber
     &=&\f{1}{c^{\f{N}{N-1}}}G(R\epsilon)\le(1-\tau\int_{B_{R\epsilon}}G^{N-1}dx\ri)
     \\\label{delt-N}
     &=&\f{1}{c^{\f{N}{N-1}}}\le(-\f{N}{\alpha_N}\log(R\epsilon)+A_0+O((R\epsilon)^N\log^{N}(R\epsilon))\ri).
     \eea
     Note that for any $T>0$, there holds
     \bna
     I_N(T)&\equiv& \int_0^T\f{t^{N-1}}{(1+t)^N}dt\\
     &=&\f{1}{1-N}\int_0^Tt^{N-1}d(1+t)^{1-N}\\
     &=&\f{1}{1-N}\le(\f{T}{1+T}\ri)^{N-1}+\int_0^T\f{t^{N-2}}{(1+t)^{N-1}}dt\\
     &=&\f{1}{1-N}\le(\f{T}{1+T}\ri)^{N-1}+I_{N-1}(T).
     \ena
     Since $I_1(T)=\log(1+T)$, we have by iteration
     $$I_N(T)=\log(1+T)-\sum_{k=1}^{N-1}\f{1}{k}\le(\f{T}{1+T}\ri)^k.$$
     Hence, by a change of variables $t=c_N(r/\epsilon)^{N(1-\beta)/(N-1)}$, we obtain
     \bea\nonumber
     \int_{B_{R\epsilon}}|\nabla\phi_\epsilon|^Ndx&=&\f{1}{\omega_{N-1}^{\f{1}{N-1}}c^{\f{N}{N-1}}}\int_0^{R\epsilon}
     \f{r^{\f{N^2}{N-1}(1-\beta)-1}}{(r^{\f{N}{N-1}(1-\beta)}+c_N^{-1}\epsilon^{\f{N}{N-1}(1-\beta)})^N}dr  \\
     \nonumber&=&\f{1}{\omega_{N-1}^{\f{1}{N-1}}c^{\f{N}{N-1}}}\f{N-1}{N(1-\beta)}\int_0^{c_NR^{\f{N}{N-1}(1-\beta)}}
     \f{t^{N-1}}{(1+t)^N}dt  \\
     \nonumber&=&\f{N-1}{\alpha_N(1-\beta)c^{\f{N}{N-1}}}\le\{\log\le(1+c_NR^{\f{N}{N-1}(1-\beta)}\ri)
     -\sum_{k=1}^{N-1}\f{1}{k}\le(\f{c_NR^{\f{N}{N-1}(1-\beta)}}{1+c_NR^{\f{N}{N-1}(1-\beta)}}\ri)^k\ri\}\\\nonumber
     &=&\f{1}{\alpha_N(1-\beta)c^{\f{N}{N-1}}}\le\{\log\f{\omega_{N-1}}{N(1-\beta)}+N(1-\beta)\log R
     \ri.\\&&\quad\quad\quad\quad\qquad\le.-(N-1)\sum_{k=1}^{N-1}\f{1}{k}
     +O(\f{1}{R^{\f{N}{N-1}(1-\beta)}})\ri\}.
     \label{delt-1-N}
     \eea
     Moreover, we require $b$ to be bounded with respect to $\epsilon$. It then follows from (\ref{cN}) that
     \be\label{delt-3N}\int_{B_{R\epsilon}}\phi_\epsilon^Ndx=O((R\epsilon)^N(\log\epsilon)^{N-1}).\ee
     Combining (\ref{delt-N})-(\ref{delt-3N}), we obtain
     \bna
     \|\phi_\epsilon\|_{1,\tau}^N&=&\f{1}{c^{\f{N}{N-1}}}\le(-\f{N}{\alpha_N}\log\epsilon+A_0-\f{N-1}{\alpha_N(1-\beta)}
     \sum_{k=1}^{N-1}\f{1}{k}
     +\f{1}{\alpha_N(1-\beta)}\log\f{\omega_{N-1}}{N(1-\beta)}\ri.\\&&\le.\qquad\quad+O(\f{1}{R^{\f{N}{N-1}(1-\beta)}})+
     O((R\epsilon)^N(\log\epsilon)^N)\ri).
     \ena
     Setting $\|\phi_\epsilon\|_{1,\tau}=1$, we have
     \be\label{c-squ-N}c^{\f{N}{N-1}}=-\f{N}{\alpha_N}\log\epsilon+A_0-\f{N-1}{\alpha_N(1-\beta)}
     \sum_{k=1}^{N-1}\f{1}{k}
     +\f{1}{\alpha_N(1-\beta)}\log\f{\omega_{N-1}}{N(1-\beta)}+O(\f{1}{R^{\f{N}{N-1}(1-\beta)}}),\ee
     which together with (\ref{cN}) leads to
     \be\label{b-2-N}b=\f{N-1}{\alpha_N(1-\beta)}\sum_{k=1}^{N-1}\f{1}{k}+O(\f{1}{R^{\f{N}{N-1}(1-\beta)}}).\ee
     Denote
     \be\label{b-e}b_\epsilon(x)=-\f{N-1}{\alpha_N(1-\beta)}\log(1+c_{N}({|x|}/{\epsilon})^{{\f{N}{N-1}(1-\beta)}})
     +b.\ee
     Then $c^{-N/(N-1)}b_\epsilon(x)=O((\log\log\epsilon^{-1})/\log\epsilon)$ uniformly in $x\in B_{R\epsilon}$, where
     $R=(\log\epsilon^{-1})^{1/(1-\beta)}$.
     We have by the Taylor formula of $(1+t)^{N/(N-1)}$ near $t=0$,
     \bea\nonumber
     \phi_\epsilon^{\f{N}{N-1}}(x)&=&c^{\f{N}{N-1}}\le(1+c^{-\f{N}{N-1}}b_\epsilon(x)\ri)^{\f{N}{N-1}}\\\nonumber
     &=&c^{\f{N}{N-1}}\le(1+\f{N}{N-1}c^{-\f{N}{N-1}}b_\epsilon(x)+\f{1}{2}\f{N}{(N-1)^2}(1+\xi)^{\f{2-N}{N-1}}
     (c^{-\f{N}{N-1}}b_\epsilon(x))^2\ri)\\\label{gt}
     &\geq&c^{\f{N}{N-1}}+\f{N}{N-1}b_\epsilon(x),
     \eea
     where $\xi$ lies between $c^{-\f{N}{N-1}}b_\epsilon(x)$ and $0$. Inserting (\ref{c-squ-N})-(\ref{b-e}) into (\ref{gt}),
     we obtain for all $x\in B_{R\epsilon}$,
     \bea\nonumber
     \alpha_N(1-\beta)\phi_\epsilon^{\f{N}{N-1}}(x)&\geq&-N(1-\beta)\log\epsilon+\alpha_N(1-\beta)A_0+\sum_{k=1}^{N-1}\f{1}{k}
     +\log\f{\omega_{N-1}}{N(1-\beta)}\\\label{BR}
     &&\qquad-N\log(1+c_{N}({|x|}/{\epsilon})^{{\f{N}{N-1}(1-\beta)}})+O(\f{1}{R^{\f{N}{N-1}(1-\beta)}}).
     \eea
     Also we have by a change of variables $t=c_Nr^{\f{N}{N-1}(1-\beta)}$ and integration by parts,
     \bea\nonumber
     \int_{B_R}\f{1}{(1+c_N|y|^{\f{N}{N-1}(1-\beta)})^N|y|^{N\beta}}dy&=&\int_0^R\f{\omega_{N-1}r^{N-1-N\beta}}
     {(1+c_Nr^{\f{N}{N-1}(1-\beta)})^N}dr\\\nonumber
     &=&\le.-\f{t^{N-2}}{(1+t)^{N-1}}\ri|_0^{c_NR^{\f{N}{N-1}(1-\beta)}}+\int_0^{c_NR^{\f{N}{N-1}(1-\beta)}}
     \f{(N-2)t^{N-3}}{(1+t)^{N-1}}dt\\\nonumber
     &=&\int_0^{c_NR^{\f{N}{N-1}(1-\beta)}}\f{1}{(1+t)^2}dt+O(\f{1}{R^{\f{N}{N-1}(1-\beta)}})\\
     &=&1+O(\f{1}{R^{\f{N}{N-1}(1-\beta)}}).\label{Inte1}
     \eea
     Combining (\ref{BR}) and (\ref{Inte1}), we obtain
     \bea\nonumber
     \int_{B_{R\epsilon}}\f{\zeta(N,\alpha_N(1-\beta)\phi_\epsilon^{\f{N}{N-1}})}
     {|x|^{N\beta}}dx&=&\int_{B_{R\epsilon}}\f{e^{\alpha_N(1-\beta)\phi_\epsilon^{\f{N}{N-1}}}}
     {|x|^{N\beta}}dx+O(c^{\f{N(N-2)}{N-1}}(R\epsilon)^{N(1-\beta)})\\\nonumber
     &\geq&\f{\omega_{N-1}}{N(1-\beta)\epsilon^{N(1-\beta)}}e^{\sum_{k=1}^{N-1}\f{1}{k}+\alpha_N(1-\beta)A_0+O(R^{-\f{N}{N-1}(1-\beta)})}\\
     \nonumber&&\times\int_{B_{R\epsilon}}\f{1}{(1+c_N(|x|/\epsilon)^{\f{N}{N-1}(1-\beta)})^N|x|^{N\beta}}dx+
     O(c^{\f{N(N-2)}{N-1}}(R\epsilon)^{N(1-\beta)})\\
     \nonumber&=&\f{\omega_{N-1}}{N(1-\beta)}e^{\sum_{k=1}^{N-1}\f{1}{k}+\alpha_N(1-\beta)A_0+O(R^{-\f{N}{N-1}(1-\beta)})}\\
     \nonumber&&\times\int_{B_R}\f{1}{(1+c_N|y|^{\f{N}{N-1}(1-\beta)})^N|y|^{N\beta}}dy+O(c^{\f{N(N-2)}{N-1}}(R\epsilon)^{N(1-\beta)})\\
     &=&\f{\omega_{N-1}}{N(1-\beta)}e^{\sum_{k=1}^{N-1}\f{1}{k}+\alpha_N(1-\beta)A_0}+O(\f{1}{R^{\f{N}{N-1}(1-\beta)}}).\label{RE}
     \eea
     Moreover,
     \bea\nonumber\int_{\mathbb{R}^N\setminus B_{R\epsilon}}
     \f{\zeta(N,\alpha_N(1-\beta)\phi_\epsilon^{\f{N}{N-1}})}{|x|^{N\beta}}dx&\geq&
     \f{\alpha_N^{N-1}(1-\beta)^{N-1}}{(N-1)! c^{\f{N}{N-1}}}\int_{\mathbb{R}^N\setminus B_{R\epsilon}}\f{G^{N}}{|x|^{N\beta}}dx\\
     &=&\f{\alpha_N^{N-1}(1-\beta)^{N-1}}{(N-1)! c^{\f{N}{N-1}}}\le(\int_{\mathbb{R}^N}\f{G^{N}}{|x|^{N\beta}}dx+o_\epsilon(1)\ri).
     \label{outer-N}\eea
     Combining (\ref{RE}), (\ref{outer-N}) and noting that $R^{-\f{N}{N-1}(1-\beta)}c^{\f{N}{N-1}}=o_\epsilon(1)$,
      we have
      $$\int_{\mathbb{R}^N}
     \f{\zeta(N,\alpha_N(1-\beta)\phi_\epsilon^{\f{N}{N-1}})}{|x|^{N\beta}}dx\geq
     \f{\omega_{N-1}e^{\sum_{k=1}^{N-1}\f{1}{k}+\alpha_N(1-\beta)A_0}}{N(1-\beta)}+
     \f{(\alpha_N(1-\beta))^{N-1}}{(N-1)! c^{\f{N}{N-1}}}\le(\int_{\mathbb{R}^N}\f{G^{N}}{|x|^{N\beta}}dx+o_\epsilon(1)\ri).$$
     Therefore we conclude (\ref{great-N}) for sufficiently small $\epsilon>0$.

     \subsection{Completion of the proof of Theorem \ref{Theorem 2}}

     Under the assumption that $c_\epsilon\ra+\infty$, there holds (\ref{u-b-N}). While it follows from (\ref{great-N}) that
     $$\Lambda_{N,\beta,\tau}>\f{1}{1-\beta}\f{\omega_{N-1}}{N}e^{\sum_{k=1}^{N-1}\f{1}{k}+\alpha_N(1-\beta)A_0}.$$
     This contradicts (\ref{u-b-N}) and implies that $c_\epsilon$ must be bounded. Then applying Theorem \ref{Estimates} to
     the equation (\ref{EL-N}), we get the desired extremal function. $\hfill\Box$

\bigskip

 {\bf Acknowledgements}. X. Li is supported by Natural Science Foundation of the Education Department of Anhui Province
 (No. KJ2016A641); Y. Yang is supported by the National Science Foundation of China (Grant Nos.11171347 and
 11471014).

\end{document}